\documentclass[11pt]{article}

\usepackage{amsmath,amsxtra,amssymb,latexsym, amscd,amsthm}
\usepackage{here}
\usepackage{graphpap}
\usepackage{graphicx}

\newtheorem{theorem}{Theorem}[section]

\newtheorem{lemma}[theorem]{Lemma}
\newtheorem{proposition}[theorem]{Proposition}
\theoremstyle{definition}
\newtheorem{definition}[theorem]{Definition}

\numberwithin{equation}{section}

\def\divv{\text{div}}
\def\<{\langle}
\def\>{\rangle}

\begin{document}


\begin{center}
\bf MATRIX COEFFICIENT IDENTIFICATION IN AN ELLIPTIC
EQUATION WITH THE CONVEX ENERGY FUNCTIONAL METHOD
\end{center}

\vspace{1cm}

\centerline {\bf Michael Hinze and Tran Nhan Tam Quyen}

{University of Hamburg, Bundesstrasse 55, 20146 Hamburg, Germany\\
Email: michael.hinze@uni-hamburg.de and quyen.tran@uni-hamburg.de}

\vspace{1cm}

{\small {\bf Abstract:} In this paper we study the inverse problem
of identifying the {\it diffusion matrix} in an elliptic PDE from measurements. The convex energy functional method with Tikhonov
regularization is applied to tackle this problem. For the discretization we use the variational discretization concept, where the PDE is discretized with piecewise linear, continuous finite elements. We show the convergence of approximations. Using a suitable source condition, we prove an error bound for discrete solutions. For the numerical solution we propose a gradient-projection algorithm and prove the strong convergence of its iterates to a solution of the identification problem. Finally, we present a numerical experiment which illustrates our theoretical results.
}


\section{Introduction}
Let $\Omega$ be an open bounded connected domain of $R^d$, $d \le 3$
with boundary $\partial\Omega$. We investigate the
problem of identifying the spatially varying diffusion matrix $Q$ in the Dirichlet
problem for the elliptic equation
\begin{align}
-\divv (Q\nabla u) &= f \mbox{ in } \Omega, \label{m1*}\\
u &=0 \mbox{ on } {\partial \Omega} \label{qmict3*}
\end{align}
from the observation $z^\delta$ of the solution $u$ in the domain $\Omega$. Here, the
function $f\in L^2(\Omega)$ is given.

In this paper we assume that $z^\delta \in
H^1_0(\Omega)$.
For related research we refer the reader to \cite{
ChanTai2003,Chavent_Kunisch2002,Cherlenyak,{Haoq},hao_quyen3,Kaltenbacher_Schoberl,kolo,wang_zou}.

Our identification problem can be considered as a generalization of
identifying the scalar diffusion coefficient $q$ in the elliptic equation
\begin{equation}\label{7/12/12:ct2}
-\divv (q\nabla u) = f \mbox{ in } \Omega \mbox{~and~} u=0 \mbox{~on~} \partial\Omega.
\end{equation}
The problem has been studied extensively in the last 30 years or so. The identification results can be found in \cite{Chicone,Know2,Ric,Vainikko-Kunisch}. Error estimates for finite element approximation solutions have been obtained, for example, in \cite{Falk,hao_quyen4,kolo,wang_zou}. A survey of numerical methods for the identification problem can be found in \cite{ChenZou,KeungZou98,Kunisch}.

Compared to the identification $q$ in (\ref{7/12/12:ct2}), the
problem of identifying the matrix $Q$ in (\ref{m1*}) has received
less attention. However, there are some contributions treating this
problem. Hoffmann and  Sprekels in \cite{Hoffmann_Sprekels} proposed a
dynamical system approach to reconstruct the matrix $Q$ in equation (\ref{m1*}).
In \cite{Rannacher_Vexler} Rannacher and
Vexler employed the finite element method and showed error estimates
for a matrix identification problem from pointwise measurements of
the state variable, provided that the sought matrix is constant and the exact data is smooth enough. 

In the present paper we adopt the convex energy functional approach of Kohn and Vogelius in \cite{Kohn_Vogelius1,Kohn_Vogelius2} to the matrix case. In fact,
for estimating the matrix $Q$ in (\ref{m1*})--(\ref{qmict3*}) from
the observation $z^\delta$ of the solution $u$, we use the {\it non-negative convex functional} (see
\S \ref{Auxiliary results})
\begin{align*}
\mathcal{J}^\delta(Q) := \int_{\Omega} Q \nabla \big( \mathcal{U}(Q)
- z^\delta\big) \cdot \nabla \big(\mathcal{U}(Q) - z^\delta \big) dx
\end{align*}
together with Tikhonov regularization and consider the
{\it strictly convex} minimization problem
$$
\min_{Q \in \mathcal{Q}_{ad}} \mathcal{J}^\delta(Q) + \rho \| Q \|^2_{{L^2(\Omega)}^{d\times d}}
$$
over the admissible set $\mathcal{Q}_{ad}$ (see \S \ref{D-I
problems}), and consider its {\it unique global} solution $Q^{\rho,\delta}$ as reconstruction. Here $\rho >0$ is the regularization parameter and $\mathcal{U}$ the non-linear
coefficient-to-solution operator.

For the discretization we use the variational discretization method introduced in \cite{Hinze} and show the convergence of approximations. Under a source condition, which is weaker than that of the existing theories in \cite{Engl_Hanke_Neubauer,EnglKuNe}, we prove an error bound for discrete regularized solutions. Finally, we employ a gradient-projection algorithm for the numerical solution of the regularized problems. The strong convergence of iterates to a solution of the identification problem is ensured without smoothness requirements on the sought matrix. Numerical results show an efficiency of our theoretical findings.

In \cite{Engl_Hanke_Neubauer,EnglKuNe} the authors investigated the convergence of Tikhonov regularized solutions via the standard output least squares method for the general non-linear ill-posed equation in Hilbert spaces. They proved some rates of convergence for this approach under a source condition and the so-called {\it small enough condition} on source elements. In the present paper, by working with a convex energy functional for our concrete identification problem, we in the proof of Theorem \ref{nu21***} are not faced with a smallness condition. Furthermore, our source condition does not require additional smoothness assumption of the sought matrix and the exact data (see \S \ref{tdht}). We also remark that such a source condition without the smallness condition was proposed in \cite{Haoq,hao_quyen1,hao_quyen2,hao_quyen3} for the {\it scalar} coefficient identification problem in elliptic PDEs and in some concrete cases the source condition was proved to satisfy if sought coefficients belong to certain smooth function spaces.

We mention that in \cite{EnglZou}, by utilizing a modified kind of adjoint, the authors for the inverse heat conduction problem introduced a source condition in the form of a variational identity without the smallness condition on source elements. The advantage of this source condition is that it does not involve the  Fr\'echet derivative of the coefficient-to-solution operator. However, the source condition requires some smoothness assumptions on the sought coefficient.

Starting with \cite{HKPS}, the authors in \cite{Grasmair,HohageWerner,WernerHohage2012} have proposed new source conditions in the form of variational inequalities. They proved some convergence rates for Tikhonov-type regularized solutions via the misfit functional method of the discrepancy for the general non-linear ill-posed equation in Banach spaces. The novelty of this theory is that the source conditions do not involve the  Fr\'echet derivative of forward operators and so avoid differentiability assumptions. Furthermore, the theory is applied to inverse problems with PDEs (see, for example, \cite{Hohage}).

Recently, by using several sets of observations and a suitable projected source condition motivated by \cite{Flemming-Hofmann} as well as certain smoothness requirements on the sought coefficient and the exact solution, the authors of \cite{Deckelnick}  derived an error bound for the finite element solutions of a standard output least squares approach to identify the diffusion matrix in \eqref{m1*}. Due to the non-linearity of the identification problem the method presented in \cite{Deckelnick} solves a non-convex minimization problem. Our approach in the present paper is different. We utilize a convex cost functional and a different source condition without smoothness assumptions. Therefore, the theory in \cite{Deckelnick} and its proof techniques are not directly comparable with our approach. Furthermore, taking the advantage of the convexity to account, we here are able to prove that iterates via a gradient-projection algorithm converge to the identified diffusion matrix.

The remaining part of this paper is organized as follows. In Section
\ref{definition} and Section \ref{Finite element method} we describe the direct and inverse problems and the finite element method which is applied to the identification problem, respectively. Convergence analysis of the finite element method is presented in Section \ref{Stability}. In Section \ref{tdht} we show
convergence rates obtained with this technique. Section \ref{iterative} is
devoted to a gradient-projection algorithm. Finally, in Section \ref{Numerical implement} we present a numerical experiment which illustrates our theoretical results.

Throughout the paper we use the
standard notion of Sobolev spaces $H^1(\Omega)$, $H^1_0(\Omega)$, $W^{k,p}(\Omega)$, etc from, for example, \cite{tro}. If not stated otherwise we write
$\int_\Omega \cdots$ instead of $\int_\Omega \cdots dx$.

\section{Problem setting and preliminaries}\label{definition}

\subsection{Notations}\label{Notation}

Let $\mathcal{S}_d$ denote the set of all symmetric $d \times
d$-matrices equipped with the inner product $M \cdot N :=
\mbox{trace} (MN)$ and the norm $$\| M \|_{\mathcal{S}_d} = (M \cdot
M)^{1/2} = \left( \sum_{i, j =1}^d m_{ij}^2 \right)^{1/2},$$ where
$M = (m_{ij})_{1\le i,j \le d}$. Let $M$ and $N$ be in
$\mathcal{S}_d$, then
$$M \preceq N $$
if and only if
$$ M \xi \cdot \xi \le N \xi \cdot \xi
~ \mbox{for all} ~ \xi \in R^d.$$
We note that if $0 \preceq M \in \mathcal{S}_d$ the root $M^{1/2}$ is well defined.

In $\mathcal{S}_d$ we introduce the convex subset
$$\mathcal{K} := \{M \in \mathcal{S}_d ~|~ \underline{q} I_d
\preceq M \preceq \overline{q} I_d\},$$
where $\underline{q}$ and $\overline{q}$ are given positive
constants and $I_d$ is the unit $d\times d$-matrix. Furthermore, let
$\xi := (\xi_1, \cdot\cdot\cdot, \xi_d)$ and $\eta :=
(\eta_1, \cdot\cdot\cdot, \eta_d)$ be two arbitrary vectors in $R^d$, we use
the notation
$$(\xi \otimes \eta)_{1\le i,j\le d} \in \mathcal{S}_d
~  \mbox{with}  ~ (\xi \otimes \eta)_{ij} := \frac{1}{2} (\xi_i \eta_j +
\xi_j \eta_i)  ~  \mbox{for all}  ~  i, j = 1, \cdots, d.$$

Finally, in the space ${L^{\infty}(\Omega)}^{d\times d}$ we use the norm
$$\|H\|_{{L^{\infty}(\Omega)}^{d\times d}} :=
\max_{1\le i,j\le d} \|h_{ij}\|_{L^{\infty}(\Omega)},$$
where $H =(h_{ij})_{1\le i,j\le d} \in {L^{\infty}(\Omega)}^{d\times d}$.

\subsection{Direct and inverse problems}\label{D-I problems}

We recall that a function $u$ in $H^1_0(\Omega)$ is said to be a weak solution of the Dirichlet problem (\ref{m1*})--(\ref{qmict3*}) if the identity
\begin{align}\label{4/6:m4}
\int_\Omega Q\nabla u \cdot \nabla v = \int_\Omega f v
\end{align}
holds for all $v\in H^1_0(\Omega)$. Assume that the matrix $Q$ belongs to the set
\begin{align}
\mathcal{Q}_{ad}:= \left\{ Q \in {L^{\infty}(\Omega)}^{d \times d}
~|~ Q(x) \in \mathcal{K}  ~  \mbox{a.e. in}  ~  \Omega\right\}.
\label{5/12/12:ct3}
\end{align}
Then, by the aid of the Poincar\'e-Friedrichs inequality in
$H^1_0(\Omega)$, there exists a positive constant $\kappa$
depending only on $\underline{q}$ and the domain $\Omega$ such that
the coercivity condition
\begin{equation}\label{coercivity}
\int_\Omega Q \nabla u \cdot \nabla u \ge \kappa
\|u\|^2_{H^1(\Omega)}
\end{equation}
holds for all $u$ in $H^1_0(\Omega)$ and $ Q \in \mathcal{Q}_{ad}$. Hence, by the
Lax-Milgram lemma, we conclude that there exists a unique solution
$u$ of (\ref{m1*})--(\ref{qmict3*}) satisfying the following estimate
\begin{align}
\left\|u\right\|_{H^1(\Omega)}\le \dfrac{1}{\kappa}
\left\|f\right\|_{L^2(\Omega)}.\label{mq5}
\end{align}

Therefore, we can define the non-linear
coefficient-to-solution operator
$$\mathcal{U} : \mathcal{Q}_{ad} \subset {L^{\infty}(\Omega)}^{d \times d}
\rightarrow H^1_0(\Omega)$$ which maps the matrix $Q \in
\mathcal{Q}_{ad} $ to the unique solution $\mathcal{U}(Q) := u$ of the problem (\ref{m1*})--(\ref{qmict3*}). Then,
the inverse problem is stated as follows:
$$\mbox{~Given~}
\overline{u} := \mathcal{U}(\overline{Q}) \in H^1_0(\Omega), \mbox{~find
a matrix~} \overline{Q} \in \mathcal{Q}_{ad} \mbox{~such that~} \eqref{4/6:m4} \mbox{~is satisfied with~} \overline{u} \mbox{~and~} \overline{Q}.$$

\subsection{Tikhonov regularization}\label{Tikhonov regularization}
According to our problem setting $\overline{u}$ is the exact solution of
(\ref{m1*})--(\ref{qmict3*}), so there exists some $\overline{Q} \in
\mathcal{Q}_{ad}$ such that $\overline{u} = \mathcal{U}(\overline{Q})$. We
assume that instead of the exact $\overline{u}$ we have only measurements $z^\delta \in H^1_0(\Omega)$ with
\begin{align}
\|z^\delta - \overline{u} \|_{H^1(\Omega)} \leq \delta \mbox{~for some~}
\delta
>0. \label{gradient-obs}
\end{align}
Our problem is to reconstruct the
matrix $\overline{Q}$ from $z^\delta$. For solving this problem we consider the
non-negative {\it convex} functional (see \S \ref{Auxiliary results})
\begin{align}
\mathcal{J}^\delta(Q) := \int_{\Omega} Q \nabla \big( \mathcal{U}(Q)
- z^\delta\big) \cdot \nabla \big(\mathcal{U}(Q) - z^\delta \big).
\label{29/6:ct8}
\end{align}
Furthermore, since the
problem is ill-posed, in this paper we shall use Tikhonov regularization
to solve it in a stable way. Namely, we consider
$$\min_{Q \in \mathcal{Q}_{ad}} \Upsilon^{\rho,\delta},
\eqno \left(\mathcal{P}^{\rho,\delta} \right)$$
where
$$\Upsilon^{\rho,\delta} := \mathcal{J}^\delta(Q) + \rho \| Q \|^2_{{L^2(\Omega)}^{d\times d}}$$
and $\rho>0$ is the regularization parameter.

In the present paper we assume that the gradient-type observation is available. Concerning this assumption, we refer the reader to \cite{hao_quyen3,Cherlenyak,Kaltenbacher_Schoberl,Chavent_Kunisch2002,ChanTai2003,kolo} and the references therein, where discussions about the interpolation of discrete measurements of the solution $\overline{u}$ which results the data $z^\delta$ satisfying \eqref{gradient-obs} are given.

\subsection{Auxiliary results}\label{Auxiliary results}

Now we summarize some properties of the coefficient-to-solution
operator. The proofs of
the following results are based on standard arguments and therefore omitted.

\begin{lemma}\label{bd21}
The coefficient-to-solution operator $\mathcal{U} :
\mathcal{Q}_{ad} \subset {L^{\infty}(\Omega)}^{d\times d}
\rightarrow H^1_0(\Omega)$ is infinitely Fr\'echet
differentiable on $\mathcal{Q}_{ad}$. For each $Q \in
\mathcal{Q}_{ad}$ and $m \ge 1$ the action of the Fr\'echet derivative $\mathcal{U}^{(m)}(Q)$ in direction $(H_1, H_2, \cdot\cdot\cdot, H_m) \in {\big({L^\infty(\Omega)}^{d \times d}\big)}^m$ denoted by $\eta := \mathcal{U}^{(m)}(Q)(H_1, H_2,
\cdot\cdot\cdot, H_m)$  is the unique weak solution in
$H^1_0(\Omega)$ to the equation
\begin{align}
\int_\Omega Q \nabla\eta \cdot \nabla v =- \sum\limits_{i=1}^m
\int_\Omega H_i\nabla \mathcal{U}^{(m-1)}(Q) \xi_i \cdot \nabla v \label{ct10}
\end{align}
for all $v\in H^1_0(\Omega)$ with $ \xi_i := (H_1,
\cdot\cdot\cdot,H_{i-1}, H_{i+1}, \cdot\cdot\cdot, H_m)$. Furthermore, the following estimate is fulfilled
$$
\|\eta\|_{H^1(\Omega)} \le\frac{m d}{\kappa^{m+1}} \|f \|_{L^2(\Omega)}
\prod_{i=1}^{m}\|H_i\|_{{L^{\infty}(\Omega)}^{d\times d}}.
$$
\end{lemma}

Now we prove the following useful results.

\begin{lemma} \label{convex}
The functional $\mathcal{J}^\delta$ defined by (\ref{29/6:ct8}) is
convex on the convex set $\mathcal{Q}_{ad}$.
\end{lemma}

\begin{proof}
From Lemma \ref{bd21} we have that $\mathcal{J}^\delta$ is infinitely differentiable. A short calculation with $\eta := \mathcal{U}'(Q)H$ gives
\begin{align*}
{\mathcal{J}^\delta}'' (Q) (H, H)
=  2\int_{\Omega} Q \nabla \eta \cdot \nabla \eta \ge 0
\end{align*}
for all $Q \in \mathcal{Q}_{ad}$ and $H \in {L^\infty(\Omega)}^{d \times d}$, which proves the lemma.
\end{proof}

\begin{lemma}[{\cite{Murat_Tartar,Tartar}}]\label{H-convergent}
Let $(Q_n)_n$ be a sequence in $\mathcal{Q}_{ad}$. Then, there exists
a subsequence, again denoted $(Q_n)_n$, and an element $Q\in
\mathcal{Q}_{ad}$ such that
\begin{quote}
$\mathcal{U}(Q_n)$ weakly converges to $\mathcal{U}(Q)$ in
$H^1_0(\Omega)$ and \\
$Q_n \nabla \mathcal{U}(Q_n)$ weakly converges to $Q\nabla
\mathcal{U}(Q)$ in ${L^2(\Omega)}^d$.
\end{quote}

The sequence $(Q_n)_n$  is then said to be H-convergent to $Q$.
\end{lemma}

The concept of H-convergence generalizes that of G-convergence introduced by Spagnolo in
\cite{Spagnolo}. Furthermore, the
H-limit of a sequence is unique.

A relationship between the H-convergence and the weak$^*$
convergence in ${L^{\infty}(\Omega)}^{d\times d}$ is given by the following lemma.

\begin{lemma} [{\cite{Murat_Tartar}}]\label{weak-convergent}
Let $(Q_n)_n$ be a sequence in $\mathcal{Q}_{ad}$. Assume that
$(Q_n)_n$  is H-convergent to $Q$ and $(Q_n)_n$ weak$^*$ converges
to $\widehat{Q}$ in ${L^\infty(\Omega)}^{d \times d}$. Then,
$Q(x) \preceq \widehat{Q}(x)$ a.e. in $\Omega$ and
\begin{align*}
\|Q\|^2_{{L^2(\Omega)}^{d \times d}} \le
\|\widehat{Q}\|^2_{{L^2(\Omega)}^{d \times d}} \le
\liminf_{n} \|Q_n\|^2_{{L^2(\Omega)}^{d \times d}}.
\end{align*}
\end{lemma}

\begin{theorem}\label{ttnghiem}
There exists a unique minimizer $Q^{\rho,\delta}$ of the
problem $\left(\mathcal{P}^{\rho,\delta} \right)$, which is called
the regularized solution of the identification problem.
\end{theorem}

\begin{proof}
Let $(Q_n)_n$ be a minimizing sequence of the problem
$(\mathcal{P}^{\rho,\delta})$, i.e., $$\lim_{n} \Upsilon^{\rho,\delta}(Q_n) = \inf_{Q\in \mathcal{Q}_{ad}}  \Upsilon^{\rho,\delta}(Q).$$ By Lemma \ref{H-convergent} and Lemma \ref{weak-convergent}, it follows that there exists
a subsequence which is not relabelled and elements $Q \in
\mathcal{Q}_{ad}$ and  $\widehat{Q} \in {L^\infty(\Omega)}^{d \times d}$ such that
\begin{quote}
$(Q_n)_n$ is H-convergent to $Q$, \\
$(Q_n)_n$ weak$^*$ converges to $\widehat{Q}$ in
${L^\infty(\Omega)}^{d \times d}$, \\
$Q(x) \preceq \widehat{Q}(x)$ a.e. in $\Omega$ and \\
$\|Q\|^2_{{L^2(\Omega)}^{d \times d}} \le
\|\widehat{Q}\|^2_{{L^2(\Omega)}^{d \times d}} \le
\liminf_{n} \|Q_n\|^2_{{L^2(\Omega)}^{d \times d}}$.
\end{quote}
We have that
\begin{align*}
\mathcal{J}^{\delta}(Q_n) &= \int_\Omega Q_n \nabla \mathcal U(Q_n) \cdot \nabla (\mathcal{U}(Q_n) - z^\delta) - \int_\Omega Q_n \nabla(\mathcal U(Q_n) - z^\delta) \cdot \nabla z^\delta\\
&= \int_\Omega f(\mathcal{U}(Q_n) - z^\delta)
- \int_\Omega Q_n \nabla \mathcal U(Q_n) \cdot \nabla z^\delta + \int_\Omega Q_n \nabla z^\delta \cdot \nabla z^\delta.
\end{align*}
And so that
\begin{align*}
\lim_n \mathcal{J}^{\delta}(Q_n)
&= \int_\Omega f(\mathcal{U}(Q) - z^\delta)
- \int_\Omega Q \nabla \mathcal U(Q) \cdot \nabla z^\delta + \int_\Omega \widehat{Q} \nabla z^\delta \cdot \nabla z^\delta\\
&= \mathcal{J}^\delta (Q) + \int_\Omega (\widehat{Q} - Q) \nabla z^\delta \cdot \nabla z^\delta \\
&\ge \mathcal{J}^\delta (Q).
\end{align*}
We therefore get
\begin{align*}
\Upsilon^{\rho,\delta}(Q)
&\le \lim_{n}
\mathcal{J}^{\delta}(Q_n) + \liminf_{n} \rho\|Q_n\|^2_{{L^2(\Omega)}^{d \times d}} \\
&= \liminf_{n} \left( \mathcal{J}^{\delta}(Q_n) + \rho\|Q_n\|^2_{{L^2(\Omega)}^{d \times d}} \right) \\
&= \inf_{Q\in \mathcal{Q}_{ad}} \mathcal{J}^{\delta}(Q) + \rho\|Q\|^2_{{L^2(\Omega)}^{d \times d}}.
\end{align*}
Since $\Upsilon^{\rho,\delta}$
is strictly convex, the minimizer is unique.
\end{proof}

\section{Discretization}\label{Finite element method}

Let $\left(\mathcal{T}_h\right)_{0<h<1}$ be a family of regular and
quasi-uniform triangulations of the domain $\overline{\Omega}$ with the mesh size $h$.
For the definition of the discretization space of the state
functions let us denote
\begin{equation*}
\mathcal{V}^1_h := \left\{v_h\in C(\overline\Omega) \cap H^1_0(\Omega)
~|~{v_h}_{|T} \in \mathcal{P}_1(T), ~~\forall
T\in \mathcal{T}_h\right\}
\end{equation*}
with $\mathcal{P}_1$ consisting all polynomial functions of degree
less than or equal to 1.
Similar to the continuous case we have the following result.
\begin{proposition}
Let $Q$ be in $\mathcal{Q}_{ad}$. Then the variational equation
\begin{align}
\int_\Omega Q\nabla u_h \cdot \nabla v_h =
\int_\Omega fv_h, \enskip\forall v_h\in
\mathcal{V}^1_h \label{10/4:ct1}
\end{align}
admits a unique solution $u_h = u_h(Q) \in \mathcal{V}^1_h$. Further, the prior estimate
\begin{align}
\left\|u_h\right\|_{H^1(\Omega)}\le \dfrac{1}{\kappa}
\left\|f\right\|_{L^2(\Omega)}, \label{18/5:ct1}
\end{align}
is satisfied.
\end{proposition}

\begin{definition}\label{discrete_solution}
The map $\mathcal{U}_h: \mathcal{Q}_{ad} \rightarrow
\mathcal{V}^1_h$ from each $Q \in  \mathcal{Q}_{ad}$
to the unique solution $u_h$ of variational equation \eqref{10/4:ct1}
is called {\it the discrete coefficient-to-solution
operator}.
\end{definition}
We note that the operator $\mathcal{U}_h$ is Fr\'echet differentiable on the set $\mathcal{Q}_{ad}$. For each $Q \in \mathcal{Q}_{ad}$ and $H \in {L^{\infty}(\Omega)}^{d\times d}$ the Fr\'echet differential $\eta_h := {\mathcal{U}_h}'(Q) H$
is an element of $\mathcal{V}_h^1$ and satisfies the equation
\begin{align}
\int_\Omega Q\nabla \eta_h \cdot \nabla v_h
&= -\int_\Omega H\nabla \mathcal{U}_h(Q) \cdot \nabla v_h  \label{ct21***}
\end{align}
for all $v_h$ in $\mathcal{V}_h^1$.

Before presenting our results we need some facts on data interpolation.

\subsection{Data interpolation}\label{Data interpolation}

It is well known that there is a usual nodal value interpolation
operator
$$I^1_h : C(\overline{\Omega}) \to \left\{v_h\in C(\overline\Omega)
~|~{v_h}_{|T} \in \mathcal{P}_1(T), ~~\forall
T\in \mathcal{T}_h\right\}$$
 such that
\begin{equation*}
I^1_h \left(H^1_0(\Omega) \cap
C(\overline{\Omega})\right) \subset \mathcal{V}^1_h.
\end{equation*}

Since $H^2(\Omega)$ is continuously embedded in $C(\overline{\Omega})$ as $d\le 3$ (see, for example, \cite{attouch}), the following result is standard in the theory of the finite element method, the
proof of which can be found, for example, in
\cite{Brenner_Scott,Ciarlet}.

\begin{lemma} \label{FEM*}
Let $\psi$ be in $H^1_0(\Omega) \cap H^2(\Omega)$. Then, we have
\begin{align*}
\left\|\psi - I^1_h \psi\right\|_{H^k(\Omega)} \le
Ch^{m-k}\left\|\psi \right\|_{H^m(\Omega)},
\end{align*}
where $0\le k< m \le 2$.
\end{lemma}

\subsection{Data mollification}\label{L2-observation}

Since the data $z^\delta$ is not smooth enough, in general we
cannot define $I^1_h$ for them, when $d \ge 2$. Instead, we use Cl\'ement's interpolation operator
$$\Pi_h: L^2(\Omega) \rightarrow \left\{v_h\in C(\overline\Omega)
~|~{v_h}_{|T} \in \mathcal{P}_1(T), ~~\forall
T\in \mathcal{T}_h\right\}$$ with
\begin{equation*}
\Pi_h \left(H^1_0(\Omega)\right) \subset \mathcal {V}^1_h
\end{equation*}
and satisfying the following convergence properties and estimates
\begin{equation}\label{23/10:ct2}
\lim_{h\to 0} \left\| \vartheta - \Pi_h \vartheta
\right\|_{H^k(\Omega)} =0 \enskip \mbox{for all} \enskip k \in \{0,
1\}
\end{equation}
and
\begin{equation}\label{23/5:ct1}
\left\| \vartheta - \Pi_h \vartheta \right\|_{H^k(\Omega)} \le
Ch^{m-k} \| \vartheta\|_{H^m(\Omega)}
\end{equation}
for $0 \le k < m \le 2$ (see \cite{Clement}, and some generalizations of which \cite{Bernardi1,Bernardi2,scott_zhang}).

\subsection{Finite element method}

Using the operator $\Pi_h$ in \S \ref{L2-observation}, we introduce
the discrete cost functional
\begin{equation}\label{29/6:ct9}
\mathcal{J}^\delta_h (Q):= \int_\Omega Q\nabla \big(\mathcal{U}_h(Q)-
\Pi_hz^{\delta}\big) \cdot  \nabla \big(\mathcal{U}_h(Q)-
\Pi_hz^{\delta}\big)
\end{equation}
with $Q \in \mathcal{Q}_{ad} $.

We note that the cost functional $\mathcal{J}^\delta_h$ contains the interpolation $\Pi_h z^\delta$ of the measurement $z^\delta$. This is different from the approaches in \cite{Deckelnick,Falk,kolo,wang_zou}. However, it is unavoidable for a numerical experiment since in general we cannot define the pointwise values of $z^\delta$ at the nodes of $\mathcal{T}_h$.

The following results are exactly obtained as in the continuous case.

\begin{lemma} \label{J-dis-convex}
For each $h>0$ the functional $\mathcal{J}^\delta_h$
defined by (\ref{29/6:ct9}) is convex on the
convex set $\mathcal{Q}_{ad}$.
\end{lemma}

We adapt a finite element version of Lemma \ref{H-convergent} and Lemma \ref{weak-convergent}.

\begin{lemma}[{\cite{Deckelnick_Hinze_2011}}]\label{Hd-convergent}
Let $(\mathcal{T}_{h_n})_n$ be a sequence of triangulations with $\lim_n h_n =0$ and $(Q_n)_n$ be a sequence in $\mathcal{Q}_{ad}$. Then there exists
a subsequence which is not relabelled and an element $Q\in
\mathcal{Q}_{ad}$ such that
\begin{quote}
$\mathcal{U}_{h_n}(Q_n)$ weakly converges to $\mathcal{U}(Q)$ in
$H^1_0(\Omega)$ and \\
$Q_n \nabla \mathcal{U}_{h_n}(Q_n)$ weakly converges to $Q\nabla
\mathcal{U}(Q)$ in ${L^2(\Omega)}^d$.
\end{quote}

The sequence $(Q_n)_n$  is then said to be Hd-convergent to $Q$.
\end{lemma}

\begin{lemma} [{\cite{Deckelnick_Hinze_2011}}]\label{d-weak-convergent}
Let $(Q_n)_n$ be a sequence in $\mathcal{Q}_{ad}$. Assume that
$(Q_n)_n$  is Hd-convergent to $Q$ and $(Q_n)_n$ weak$^*$ converges
to $\widehat{Q}$ in ${L^\infty(\Omega)}^{d \times d}$. Then,
$Q(x) \preceq \widehat{Q}(x)$ a.e. in $\Omega$ and
\begin{align*}
\|Q\|^2_{{L^2(\Omega)}^{d \times d}} \le
\|\widehat{Q}\|^2_{{L^2(\Omega)}^{d \times d}} \le
\liminf_{n} \|Q_n\|^2_{{L^2(\Omega)}^{d \times d}}.
\end{align*}
\end{lemma}

\begin{lemma} \label{solution2}
Let $$\Upsilon^{\rho,\delta}_h (Q):=
\mathcal{J}^\delta_h(Q) + \rho  \left \| Q
\right \|^2_{{L^2(\Omega)}^{d\times d}}.$$
There exists a unique minimizer $Q^{\rho,\delta}_h$ of the strictly convex minimization problem
$$
\min_{Q\in \mathcal{Q}_{ad}} \Upsilon^{\rho,\delta}_h (Q). \eqno \left(
\mathcal{P}^{\rho,\delta}_h \right)
$$
\end{lemma}

Now we consider the orthogonal projection
$P_{\mathcal{K}} : \mathcal{S}_d \to \mathcal{K}$ characterised by
$$(A - P_{\mathcal{K}}(A)) (B - P_{\mathcal{K}}(A)) \le 0$$
for all $A \in \mathcal{S}_d$ and $B \in \mathcal{K}$.

\begin{lemma}\label{Projection}
Let $Q^{\rho,\delta}_h \in \mathcal{Q}_{ad}$. Then
$Q^{\rho,\delta}_h$ is the unique solution  of the problem
$\left(\mathcal{P}^{\rho,\delta}_h \right)$ if and only if the equation
\begin{align*}
Q^{\rho,\delta}_h(x) = P_{\mathcal{K}}
\left(\frac{1}{2\rho} \left( \nabla \mathcal{U}_h(Q^{\rho,\delta}_h)(x)
\otimes \nabla \mathcal{U}_h(Q^{\rho,\delta}_h)(x)
- \nabla  \Pi_hz^{\delta} (x) \otimes \nabla  \Pi_hz^{\delta} (x) \right)\right)
\end{align*}
holds for a.e. in $\Omega$.
\end{lemma}
\begin{proof}
Since the problem $\left(
\mathcal{P}^{\rho,\delta}_h \right)$ is strictly convex,
an element $Q^{\rho,\delta}_h \in \mathcal{Q}_{ad}$ is
the unique solution  of $\left(
\mathcal{P}^{\rho,\delta}_h \right)$ if and only if the inequality
\begin{align}\label{17/10/14:ct5}
{\mathcal{J}^\delta_h}'\big(Q^{\rho,\delta}_h\big) \big(Q - Q^{\rho,\delta}_h\big) + 2\rho \int_\Omega Q^{\rho,\delta}_h
\cdot \big(Q - Q^{\rho,\delta}_h\big) \ge 0
\end{align}
is satisfied for all $Q \in \mathcal{Q}_{ad}$. 

By (\ref{29/6:ct9}) and (\ref{ct21***}), we have that
\begin{align} \label{17/10/14:ct6}
{\mathcal{J}^\delta_h}'\big(Q^{\rho,\delta}_h\big) \big(Q - Q^{\rho,\delta}_h\big)
&= \int_\Omega \big(Q - Q^{\rho,\delta}_h\big)\nabla \big(\mathcal{U}_h \big(Q^{\rho,\delta}_h\big) -
\Pi_hz^{\delta}\big) \cdot  \nabla \big(\mathcal{U}_h \big(Q^{\rho,\delta}_h\big) - \Pi_hz^{\delta}\big) \nonumber \\
&~\quad + 2\int_\Omega Q^{\rho,\delta}_h \nabla {\mathcal{U}_h}' \big(Q^{\rho,\delta}_h\big)\big(Q - Q^{\rho,\delta}_h\big) \cdot  \nabla \big(\mathcal{U}_h \big(Q^{\rho,\delta}_h\big) - \Pi_hz^{\delta}\big) \nonumber\\
&= \int_\Omega \big(Q - Q^{\rho,\delta}_h\big)\nabla \big(\mathcal{U}_h \big(Q^{\rho,\delta}_h\big) -
\Pi_hz^{\delta}\big) \cdot  \nabla \big(\mathcal{U}_h \big(Q^{\rho,\delta}_h\big) - \Pi_hz^{\delta}\big) \nonumber\\
&~\quad - 2\int_\Omega \big(Q - Q^{\rho,\delta}_h\big) \nabla \mathcal{U}_h \big(Q^{\rho,\delta}_h\big) \cdot  \nabla \big(\mathcal{U}_h \big(Q^{\rho,\delta}_h\big) - \Pi_hz^{\delta}\big) \nonumber\\
&= - \int_\Omega \big(Q - Q^{\rho,\delta}_h\big) \left( \nabla \mathcal{U}_h \big(Q^{\rho,\delta}_h\big) \cdot  \nabla \mathcal{U}_h \big(Q^{\rho,\delta}_h\big) - \nabla \Pi_hz^{\delta} \cdot \nabla \Pi_hz^{\delta}\right) \nonumber\\
&= - \int_\Omega \left( \nabla \mathcal{U}_h(Q^{\rho,\delta}_h)
\otimes \nabla \mathcal{U}_h(Q^{\rho,\delta}_h)
- \nabla  \Pi_hz^{\delta} \otimes \nabla  \Pi_hz^{\delta} \right)
\cdot \big(Q - Q^{\rho,\delta}_h\big).
\end{align}
It follows from (\ref{17/10/14:ct5}) and (\ref{17/10/14:ct6}) that
\begin{align*}
\int_\Omega \left(\frac{1}{2\rho} \left( \nabla \mathcal{U}_h \big(Q^{\rho,\delta}_h\big) \otimes \nabla
\mathcal{U}_h \big(Q^{\rho,\delta}_h\big)
- \nabla  \Pi_hz^{\delta} \otimes \nabla  \Pi_hz^{\delta}\right) - Q^{\rho,\delta}_h\right)
\cdot \big(Q - Q^{\rho,\delta}_h\big) \le 0
\end{align*}
for all $Q \in \mathcal{Q}_{ad}$. Then a localization argument infers
\begin{align*}
\left( \frac{1}{2\rho} \left( \nabla \mathcal{U}_h \big(Q^{\rho,\delta}_h\big) (x) \otimes \nabla \mathcal{U}_h \big(Q^{\rho,\delta}_h\big)(x)
- \nabla  \Pi_hz^{\delta}(x) \otimes \nabla  \Pi_hz^{\delta}(x) \right)
- Q^{\rho,\delta}_h (x) \right) \cdot \\
\cdot \big(M - Q^{\rho,\delta}_h (x)\big)
\le 0 \quad \quad \quad \quad \quad\enskip
\end{align*}
for all $M \in \mathcal{K}$. The proof is completed.
\end{proof}

{\bf Remark.} Since $\mathcal{U}_h \big(Q^{\rho,\delta}_h\big)$ and $\Pi_h z^\delta$ are both in $\mathcal{V}^1_h$, the assertion of Lemma \ref{Projection} shows that the solution of
$\left( \mathcal{P}^{\rho,\delta}_h \right)$
is a piecewise constant matrix over $\mathcal{T}_h$, so that it belongs to the set
$\mathcal{Q}_{ad} \cap \mathcal{V}_h$, where
\begin{align} \label{18/10/14:ct1}
\mathcal{V}_h:= \Big\{ M := (m_{ij})_{1\le i,j\le d} \in & {L^{\infty}(\Omega)}^{d \times d}
~ \Big|~ M(x) \in \mathcal{S}_d  \mbox{~a.e. in~} \Omega \mbox{~and~} \notag\\
& {m_{ij}}_{|T} = \mbox{const} \mbox{~for all~} i,j \mbox{~with~}
1\le i,j\le d \mbox{~and~} T\in \mathcal{T}_h\Big\}.
\end{align}
Taking this into account, a discretization of the admissible set
$\mathcal{Q}_{ad}$ can be avoided. Furthermore,
we note that $\mathcal{Q}_{ad} \cap \mathcal{V}_h$
is a non-empty, convex, bounded and closed set in the
${L^2(\Omega)}^{d \times d}$-norm in the finite
dimensional space $\mathcal{V}_h$.

In what follows $C$ is a generic positive constant which is independent of the mesh size $h$ of $\mathcal{T}_h$, the noise level $\delta$ and the regularization parameter $\rho$.

\section{Convergence} \label{Stability}

In this section we analyze the convergence of Tikhonov regularization. To this end, we introduce
\begin{align}\label{19-6-15ct1}
\sigma_h(Q) := \| \mathcal{U}(Q) - \mathcal{U}_h(Q) \|_{H^1(\Omega)}
\end{align}
and
\begin{align}\label{19-6-15ct1*}
\gamma_h(\varphi) := \| \varphi - \Pi_h\varphi \|_{H^1(\Omega)},
\end{align}
where $Q \in \mathcal{Q}_{ad}$ and $\varphi \in H^1_0(\Omega)$. We note that
$$\lim_{h \to 0} \sigma_h(Q)  = 0 \mbox{~ and ~} \lim_{h \to 0} \gamma_h(\varphi) = 0.$$

\begin{theorem}\label{odinh1}
Let $\left(\mathcal{T}_{h_n}\right)_n$ be a sequence of triangulations with
$\lim_n h_n = 0$. Assume that
$\big( Q^{\rho, \delta}_{h_n} \big)_n$ is the sequence of unique minimizers of
$\big(\mathcal{P}^{\rho,\delta}_{h_n} \big)$.
Then $\big( Q^{\rho, \delta}_{h_n} \big)_n$ converges to
minimizer $Q^{\rho, \delta}$ of
$\big( \mathcal{P}^{\rho,\delta} \big)$
in the ${L^2(\Omega)}^{d\times d}$-norm.
\end{theorem}

\begin{proof}
For the sake of the notation we denote by
$$Q_n := Q^{\rho, \delta}_{h_n}.$$
In view of Lemma \ref{Hd-convergent} and Lemma \ref{d-weak-convergent} there exists a subsequence, again denoted $(Q_n)_n$ and elements $Q^{\rho,\delta} \in \mathcal{Q}_{ad}$, $\widehat{Q} \in {L^\infty(\Omega)}^{d \times d}$ such that $(Q_n)_n$ is Hd-convergent to $Q^{\rho,\delta}$, $(Q_n)_n$ weak$^*$ converges to $\widehat{Q}$ in
${L^\infty(\Omega)}^{d \times d}$, $Q^{\rho,\delta}(x) \preceq \widehat{Q}(x)$ a.e. in $\Omega$ and $\|Q^{\rho,\delta}\|^2_{{L^2(\Omega)}^{d \times d}} \le \|\widehat{Q}\|^2_{{L^2(\Omega)}^{d \times d}} \le \liminf_{n} \|Q_n\|^2_{{L^2(\Omega)}^{d \times d}}$.

First we show that
\begin{align}\label{19-6-15ct2}
\lim_n \mathcal{J}_{h_n}^{\delta}\big(Q_n\big) \ge \mathcal{J}^\delta\big(Q^{\rho,\delta}\big).
\end{align}
Indeed, we write
\begin{align}\label{19-6-15ct3}
\mathcal{J}_{h_n}^{\delta}\big(Q_n\big)
&= \int_\Omega Q_n \nabla \mathcal{U}_{h_n}(Q_n) \cdot  \nabla \mathcal{U}_{h_n}(Q_n) - 2\int_\Omega Q_n \nabla \mathcal{U}_{h_n}(Q_n)\cdot  \nabla  \Pi_{h_n}z^{\delta}
+ \int_\Omega Q_n \nabla \Pi_{h_n}z^{\delta} \cdot \nabla \Pi_{h_n}z^{\delta}.
\end{align}
We have that
\begin{align}\label{22-6-15ct1}
\int_\Omega Q_n \nabla \mathcal{U}_{h_n}(Q_n) \cdot  \nabla \mathcal{U}_{h_n}(Q_n) &= \int_\Omega f \mathcal{U}_{h_n}(Q_n)\notag\\
& \rightarrow \int_\Omega f \mathcal{U}(Q^{\rho,\delta})
\end{align}
and, by \eqref{23/10:ct2},
\begin{align}\label{22-6-15ct2}
\int_\Omega Q_n \nabla \mathcal{U}_{h_n}(Q_n)\cdot  \nabla  \Pi_{h_n}z^{\delta} &= \int_\Omega Q_n \nabla \mathcal{U}_{h_n}(Q_n)\cdot  \nabla  z^{\delta} + \int_\Omega Q_n \nabla \mathcal{U}_{h_n}(Q_n)\cdot  \nabla  \big(\Pi_{h_n}z^{\delta} - z^\delta\big)\notag\\
& \rightarrow \int_\Omega Q^{\rho,\delta}\nabla \mathcal{U}(Q^{\rho,\delta})\cdot  \nabla  z^{\delta}
\end{align}
and
\begin{align}\label{22-6-15ct3}
\int_\Omega Q_n \nabla \Pi_{h_n}z^{\delta} \cdot \nabla \Pi_{h_n}z^{\delta} &= \int_\Omega Q_n \nabla z^{\delta} \cdot \nabla z^{\delta} + \int_\Omega Q_n \nabla (\Pi_{h_n}z^{\delta} - z^\delta) \cdot \nabla \Pi_{h_n}z^{\delta} \notag\\
&~\quad +  \int_\Omega Q_n \nabla (\Pi_{h_n}z^{\delta} - z^\delta) \cdot \nabla z^{\delta} \notag\\
&\rightarrow \int_\Omega \widehat{Q} \nabla z^{\delta} \cdot \nabla z^{\delta} \notag\\
&= \int_\Omega Q^{\rho,\delta} \nabla z^{\delta} \cdot \nabla z^{\delta} + \int_\Omega (\widehat{Q} - Q^{\rho,\delta})\nabla z^{\delta} \cdot \nabla z^{\delta} \notag\\
&\ge \int_\Omega Q^{\rho,\delta} \nabla z^{\delta} \cdot \nabla z^{\delta}.
\end{align}
By \eqref{19-6-15ct3}--\eqref{22-6-15ct3}, we arrive at \eqref{19-6-15ct2}. Furthermore, in view of \eqref{19-6-15ct1} and \eqref{23/10:ct2}, for all $Q \in \mathcal{Q}_{ad}$ we also get
\begin{align*}
\lim_n \mathcal{J}_{h_n}^{\delta}(Q) = \mathcal{J}^\delta(Q).
\end{align*}
Hence it follows that for all $Q \in \mathcal{Q}_{ad}$
\begin{align}\label{22-6-15ct4}
\mathcal{J}^\delta\big( Q^{\rho,\delta} \big) + \rho \big\|
Q^{\rho,\delta}  \big\|^2_{{L^2(\Omega)}^{d\times d}} &\le \lim_n \mathcal{J}_{h_n}^{\delta}(Q_n) + \liminf_n \rho \big\| Q_n \big\|^2_{{L^2(\Omega)}^{d\times d}} \notag\\
&= \liminf_n \left( \mathcal{J}_{h_n}^{\delta}(Q_n)
+ \rho \big\| Q_n\big\|^2_{{L^2(\Omega)}^{d\times d}} \right) \nonumber \\
&\le \limsup_n \left( \mathcal{J}_{h_n}^{\delta}(Q_n)
+ \rho \big\| Q_n
\big\|^2_{{L^2(\Omega)}^{d\times d}} \right) \nonumber \\
& \le \limsup_n \left( \mathcal{J}_{h_n}^{\delta} (Q) + \rho \big\| Q
 \big\|^2_{{L^2(\Omega)}^{d\times d}} \right) \nonumber\\
& = \lim_n \left( \mathcal{J}_{h_n}^{\delta} (Q) + \rho \big\| Q
 \big\|^2_{{L^2(\Omega)}^{d\times d}} \right) \nonumber\\
&= \mathcal{J}^\delta (Q) + \rho \big\| Q
\big\|^2_{{L^2(\Omega)}^{d\times d}}.
\end{align}
Thus, $Q^{\rho,\delta}$ is a unique solution to
$\big(\mathcal{P}^{\rho,\delta}\big)$. It remains to show that $( Q_n)_n$ converges to
$Q^{\rho, \delta}$ in the ${L^2(\Omega)}^{d\times d}$-norm. To this end, we rewrite
\begin{align*}
\rho \big\| Q^{\rho,\delta} - Q_n \big\|^2_{{L^2(\Omega)}^{d\times d}}
&= \rho \big\| Q^{\rho,\delta} \|^2_{{L^2(\Omega)}^{d\times d}} + \rho \big\| Q_n \big\|^2_{{L^2(\Omega)}^{d\times d}}  -2 \rho \big\langle Q^{\rho,\delta} , Q_n \big\rangle_{{L^2(\Omega)}^{d\times d}} \notag\\
&= \rho \big\| Q^{\rho,\delta} \|^2_{{L^2(\Omega)}^{d\times d}}  + \left( \mathcal{J}_{h_n}^{\delta}(Q_n) + \rho \big\| Q_n \big\|^2_{{L^2(\Omega)}^{d\times d}} \right)  \notag\\
&~\quad -2 \rho \big\langle Q^{\rho,\delta} , Q_n \big\rangle_{{L^2(\Omega)}^{d\times d}} - \mathcal{J}_{h_n}^{\delta}(Q_n).
\end{align*}
By \eqref{22-6-15ct4}, we have that
$$\limsup_n \left( \mathcal{J}_{h_n}^{\delta}(Q_n)
+ \rho \big\| Q_n \big\|^2_{{L^2(\Omega)}^{d\times d}} \right) = \mathcal{J}^\delta\big( Q^{\rho,\delta} \big) + \rho \big\|
Q^{\rho,\delta}  \big\|^2_{{L^2(\Omega)}^{d\times d}}.$$
Therefore, by \eqref{19-6-15ct2} and the fact that $(Q_n)_n$ weak$^*$ converges to $\widehat{Q}$ in ${L^\infty(\Omega)}^{d \times d}$ with $Q^{\rho,\delta}(x) \preceq \widehat{Q}(x)$ a.e. in $\Omega$, we deduce that
\begin{align*}
\rho \lim_n \big\| Q^{\rho,\delta} - Q_n \big\|^2_{{L^2(\Omega)}^{d\times d}}
&\le \rho \big\| Q^{\rho,\delta} \|^2_{{L^2(\Omega)}^{d\times d}}  + \left( \mathcal{J}^\delta\big( Q^{\rho,\delta} \big) + \rho \big\|
Q^{\rho,\delta}  \big\|^2_{{L^2(\Omega)}^{d\times d}} \right)  \notag\\
&~\quad -2 \rho \big\langle Q^{\rho,\delta} , \widehat{Q} \big\rangle_{{L^2(\Omega)}^{d\times d}} - \mathcal{J}^\delta\big( Q^{\rho,\delta} \big)\\
&\le 2\rho \big\| Q^{\rho,\delta} \|^2_{{L^2(\Omega)}^{d\times d}} -2 \rho \big\langle Q^{\rho,\delta} , Q^{\rho,\delta} \big\rangle_{{L^2(\Omega)}^{d\times d}} \\
&= 0.
\end{align*}
The proof is completed.
\end{proof}

Next we show convergence of
discrete regularized solutions to identification problem. Before presenting our result we introduce the notion of the minimum norm solution of the
identification problem.

\begin{lemma} \label{nx45}
The set
$$
\mathcal{I}_{\mathcal{Q}_{ad}}(\overline{u}) := \{ Q\in
\mathcal{Q}_{ad}~|~\mathcal{U}(Q)=\overline{u}\}
$$
is non-empty, convex, bounded and closed in the ${L^2(\Omega)}^{d\times d}$-norm. Hence
there is a unique minimizer $Q^\dag$ of the
problem
\begin{align*}
\min_{ Q \in \mathcal{I}_{\mathcal{Q}_{ad}}(\overline{u})}\left\| Q\right\|^2_{{L^2(\Omega)}^{d\times d}}
\end{align*}
which is called by the minimum norm solution of the
identification problem.
\end{lemma}

\begin{theorem}\label{convergence1}
Let $\left(\mathcal{T}_{h_n}\right)_n$ be a sequence of triangulations  with
mesh sizes $\left(h_n\right)_n$. Let $\left(\delta_n\right)_n$
and $\left(\rho_n\right)_n$ be any positive sequences such that
$$\rho_n\rightarrow 0, ~\frac{\delta_n^2}{\rho_n}
\rightarrow 0, ~\frac{\sigma^2_{h_n}(Q^\dag)}{\rho_n}
\rightarrow 0 \mbox{~and~} \frac{\gamma^2_{h_n} \big(\mathcal{U}(Q^\dag)\big)}{\rho_n}
\rightarrow 0$$
as $n\rightarrow\infty$. Moreover, assume that $\left(z^{\delta_n}\right)_n $ is a sequence satisfying $\left \| \mathcal{U}(Q^\dag) - z^{\delta_n}
\right \|_{H^1_0(\Omega)} \le \delta_n$ and $\big( Q^{\rho_n, \delta_n}_{h_n} \big)_n$ is the sequence of unique minimizers of $\big( \mathcal{P}^{\rho_n,\delta_n}_{h_n} \big)$. Then $\big(Q^{\rho_n, \delta_n}_{h_n}\big)_n$
converges to $Q^\dag$ in the ${L^2(\Omega)}^{d\times d}$-norm as $n\to \infty$.
\end{theorem}

\begin{proof}
Denoting
$$Q_n := Q^{\rho_n, \delta_n}_{h_n}$$
and due to the definition of $Q_n$, we get
\begin{align*}
\mathcal{J}^{\delta_n}_{h_n} \big(Q_n\big) + \rho_n
\big \| Q_n  \big\|^2_{{L^2(\Omega)}^{d\times d}} &\le
\mathcal{J}^{\delta_n}_{h_n} (Q^\dag) + \rho_n \big\| Q^\dag
\big\|^2_{{L^2(\Omega)}^{d\times d}}.
\end{align*}
We have that
\begin{align*}
\mathcal{J}^{\delta_n}_{h_n}& (Q^\dag) \\
&= \int_\Omega Q^\dag \nabla \big(\mathcal{U}_{h_n}(Q^\dag)- \Pi_{h_n} z^{\delta_n}\big) \cdot \big(\mathcal{U}_{h_n}(Q^\dag)- \Pi_{h_n} z^{\delta_n}\big) \\
&\le \overline{q} \big\| \mathcal{U}_{h_n}(Q^\dag)- \Pi_{h_n}z^{\delta_n}\big\|^2_{H^1(\Omega)} \\
&= \overline{q} \big\| \mathcal{U}_{h_n}(Q^\dag) - \mathcal{U}(Q^\dag) + \Pi_{h_n} \big(\mathcal{U}(Q^\dag) - z^{\delta_n}\big) + \mathcal{U}(Q^\dag) - \Pi_{h_n} \mathcal{U}(Q^\dag)\big\|^2_{H^1(\Omega)} \\
&\le 3\overline{q} \left( \big\| \mathcal{U}_{h_n}(Q^\dag) - \mathcal{U}(Q^\dag) \big\|^2_{H^1(\Omega)} + \big\| \Pi_{h_n} \big(\mathcal{U}(Q^\dag) - z^{\delta_n} \big)\big\|^2_{H^1(\Omega)} + \big\| \mathcal{U}(Q^\dag) - \Pi_{h_n} \mathcal{U}(Q^\dag)\big\|^2_{H^1(\Omega)} \right) \\
& \le 3C\overline{q} \left( \sigma^2_{h_n}(Q^\dag) + \delta_n^2 + \gamma^2_{h_n}\big(\mathcal{U}(Q^\dag)\big)\right),
\end{align*}
where $C$ is the positive constant defined by
$$C := \max \big(1, \| \Pi_{h_n} \|_{\mathcal{L}(H^1(\Omega), H^1(\Omega))} \big).$$
So that
\begin{align}
\mathcal{J}^{\delta_n}_{h_n} \big(Q_n\big) + \rho_n
\big \| Q_n  \big\|^2_{{L^2(\Omega)}^{d\times d}} &\le
3C\overline{q} \big( \sigma^2_{h_n}(Q^\dag) + \delta_n^2 + \gamma^2_{h_n}\big(\mathcal{U}(Q^\dag)\big) + \rho_n \big\| Q^\dag
\big\|^2_{{L^2(\Omega)}^{d\times d}}.\label{odinh22}
\end{align}
This implies
\begin{align}
\limsup_n\big \| Q_n \big\|^2_{{L^2(\Omega)}^{d\times d}} &\le \limsup_n
\left( 3C\overline{q} \frac{\sigma^2_{h_n}(Q^\dag) + \delta_n^2 + \gamma^2_{h_n}\big(\mathcal{U}(Q^\dag)\big)}{\rho_n} +
\big\| Q^\dag  \big\|^2_{{L^2(\Omega)}^{d\times d}} \right)\notag\\
&= \big\| Q^\dag
\big\|^2_{{L^2(\Omega)}^{d\times d}}.\label{odinh22*}
\end{align}
By Lemma \ref{Hd-convergent} and Lemma \ref{d-weak-convergent} there exists a subsequence which is not relabelled and elements $\Theta \in
\mathcal{Q}_{ad}$, $\widehat{Q} \in {L^\infty(\Omega)}^{d \times d}$ such that $(Q_n)_n$ is Hd-convergent to $\Theta$, $(Q_n)_n$ weak$^*$ converges to $\widehat{Q}$ in
${L^\infty(\Omega)}^{d \times d}$, $\Theta(x) \preceq \widehat{Q}(x)$ a.e. in $\Omega$ and
\begin{align}\label{23-6-15ct1}
\|\Theta\|^2_{{L^2(\Omega)}^{d \times d}} \le \|\widehat{Q}\|^2_{{L^2(\Omega)}^{d \times d}} \le \liminf_{n} \|Q_n\|^2_{{L^2(\Omega)}^{d \times d}}.
\end{align}
Moreover, the proof of Theorem \ref{odinh1} includes an argument which can be used to show that
\begin{align*}
\lim_n \mathcal{J}_{h_n}^{\delta_n}\big(Q_n\big) \ge \int_{\Omega} \Theta \nabla \big(
\mathcal{U}(\Theta) - \mathcal{U}(Q^\dag)\big) \cdot
\nabla \big( \mathcal{U}(\Theta) - \mathcal{U}(Q^\dag)\big).
\end{align*}
Then by (\ref{coercivity}) and (\ref{odinh22}), we arrive at
\begin{align*}
\kappa \big\| \mathcal{U}(\Theta) - \mathcal{U}(Q^\dag)
\big\|^2_{H^1(\Omega)} \le \lim_n \mathcal{J}_{h_n}^{\delta_n}\big(Q_n\big) = 0.
\end{align*}
Therefore, $\Theta \in\mathcal{I}_{\mathcal{Q}_{ad}}(\overline{u}).$ Furthermore, by \eqref{23-6-15ct1} and \eqref{odinh22*} and the uniqueness of the minimum norm solution, we obtain that $\Theta = Q^\dag$. Finally, by \eqref{odinh22*}, the fact that $(Q_n)_n$ weak$^*$ converges to $\widehat{Q}$ in ${L^\infty(\Omega)}^{d \times d}$ and $Q^\dag(x) \preceq \widehat{Q}(x)$ a.e. in $\Omega$, we infer that
\begin{align*}
\limsup_n\big \| Q_n - Q^\dag \big\|^2_{{L^2(\Omega)}^{d\times d}} &= \limsup_n \left( \big \| Q_n \big\|^2_{{L^2(\Omega)}^{d\times d}} + \big \| Q^\dag \big\|^2_{{L^2(\Omega)}^{d\times d}} - 2 \big\langle Q_n, Q^\dag \big\rangle_{{L^2(\Omega)}^{d\times d}}\right) \\
&\le 2\big \| Q^\dag \big\|^2_{{L^2(\Omega)}^{d\times d}} - 2 \big\langle \widehat{Q}, Q^\dag \big\rangle_{{L^2(\Omega)}^{d\times d}}\\
&\le 2\big \| Q^\dag \big\|^2_{{L^2(\Omega)}^{d\times d}} - 2 \big\langle Q^\dag, Q^\dag \big\rangle_{{L^2(\Omega)}^{d\times d}} \\
&=0.
\end{align*}
The proof is completed.
\end{proof}

\section{Convergence rates}\label{tdht}

Now we state the result on convergence rates for Tikhonov
regularization of our identification problem. Before presenting the result we recall some notions.

Any $\Psi \in L^\infty(\Omega)$ can be considered as an element in ${L^\infty(\Omega)}^*$ by
\begin{align}
\left\<\Psi, \psi\right\>_{\left({L^\infty(\Omega)}^*,
L^\infty(\Omega)\right)} := \int_\Omega \Psi \psi  \label{inf1}
\end{align}
for all $\psi \mbox{~in~} L^\infty(\Omega)$ with
$\| \Psi \|_{{L^{\infty}(\Omega)}^*}\le |\Omega| \| \Psi \|_{L^\infty(\Omega)}$.

According to Lemma \ref{bd21}, for each $Q \in \mathcal{Q}_{ad}$ the mapping
$$
\mathcal{U}'(Q): {L^{\infty}(\Omega)}^{d \times d} \rightarrow H^1_0(\Omega)
$$
is a continuous operator with the dual
$$
{\mathcal{U}'(Q)}^*:{H^1_0(\Omega)}^* \rightarrow
{\left({L^{\infty}(\Omega)}^{d \times d}\right)}^*.
$$
Let $w^* \in {H^1_0(\Omega)}^*$ be arbitrary but fixed. We consider the Dirichlet problem
\begin{align}\label{25-9-15ct1}
-\divv (Q^\dag\nabla w) = w^* \mbox{ in } \Omega \mbox{~and~} w=0 \mbox{~on~} \partial\Omega
\end{align}
which has a unique weak solution  $w \in H^1_0(\Omega)$. Then for all $H \in {L^{\infty}(\Omega)}^{d \times d}$ we have
\begin{align} \label{29-6-15ct1}
\big \langle {\mathcal{U}'(Q^\dag)}^* w^*, H \big \rangle_{\big({{L^{\infty}(\Omega)}^{d \times d}}^*, {L^{\infty}(\Omega)}^{d \times d}\big)}
&= \langle w^*, \mathcal{U}'(Q^\dag)H \rangle_{\big({H^1_0(\Omega)}^*, H^1_0(\Omega)\big)} \notag\\
&= \int_\Omega Q^\dag \nabla w \cdot \nabla \mathcal{U}'(Q^\dag)H.
\end{align}

\begin{theorem}\label{nu21***}
Assume that there is a functional $w^* \in {H^1_0 (\Omega)}^*$ such that
\begin{align}\label{moi14***}
{\mathcal{U}'(Q^\dag)}^*w^* = Q^\dag.
\end{align}
Then
\begin{align} \label{29-6-15ct3}
\dfrac{\kappa}{4}\big\| \mathcal{U}_h(Q_h) - \mathcal{U} (Q^\dag) \big\|^2_{H^1(\Omega)} + \rho\big\|Q_h - Q^\dag \big\|^2_{{L^2(\Omega)}^{d\times d}} =\mathcal{O} \left(\delta^2 + \sigma_h^2(Q^\dag) + \gamma_h^2 \big( \mathcal{U}(Q^\dag) \big) + \gamma^2_h (w) + \rho^2 \right),
\end{align}
where $Q_h := Q^{\rho,\delta}_h$ is the unique solution of $\big(
\mathcal{P}^{\rho,\delta}_h \big)$.
\end{theorem}

We remark that in case $\overline{u} ,w \in H^2(\Omega)$ with $w$ from \eqref{25-9-15ct1}, by the C\'ea's lemma and \eqref{23/5:ct1}, we infer that $\sigma_h(Q^\dag) \le Ch$, $\gamma_h\big( \mathcal{U}(Q^\dag) \big) \le Ch$ and $\gamma_h(w) \le Ch$. Therefore, the convergence rate
$$\dfrac{\kappa}{4}\big\| \mathcal{U}_h(Q_h) - \mathcal{U} (Q^\dag) \big\|^2_{H^1(\Omega)} + \rho\big\|Q_h - Q^\dag \big\|^2_{{L^2(\Omega)}^{d\times d}} =\mathcal{O} \left(\delta^2 + h^2 + \rho^2 \right)$$
is obtained.

By \eqref{ct10}, \eqref{inf1} and \eqref{29-6-15ct1}, the source condition \eqref{moi14***} is satisfied if there exists a functional $w \in H^1_0(\Omega)$ such that for all $H \in {L^{\infty}(\Omega)}^{d \times d}$ the equation
$$\int_\Omega H \cdot Q^\dag = - \int_\Omega H \nabla \mathcal{U}(Q^\dag) \cdot \nabla w$$
holds. However, as we can see in \eqref{28-9-15ct1} below, the convergence rate \eqref{29-6-15ct3} is obtained under the weaker condition that there exists a functional $w \in H^1_0(\Omega)$ such that
\begin{align}\label{28-9-15ct2}
\int_\Omega (Q^\dag - Q_h) \cdot Q^\dag \le \left| \int_\Omega (Q^\dag - Q_h) \nabla \mathcal{U}(Q^\dag) \cdot \nabla w\right|.
\end{align}

\begin{lemma}
If there exists a functional $w \in H^1_0(\Omega)$ such that
\begin{align}\label{28-9-15ct3}
Q^\dag(x) = P_{\mathcal{K}}\left(\nabla \mathcal{U}(Q^\dag) (x) \otimes \nabla w(x)\right) \mbox{~a.e in~} \Omega,
\end{align}
then the condition \eqref{28-9-15ct2} holds. Thus the convergence rate \eqref{29-6-15ct3} is obtained.
\end{lemma}

We note that \eqref{28-9-15ct3} is the projected source condition introduced in \cite{Deckelnick}. However, we here do not require any of the smoothness of the sought matrix and the exact data.

\begin{proof}
We have
\begin{align*}
\int_\Omega (Q^\dag - Q_h) \cdot Q^\dag &= - \int_\Omega P_{\mathcal{K}}\left(\nabla \mathcal{U}(Q^\dag) \otimes \nabla w\right) \cdot \left(Q_h - P_{\mathcal{K}}\left(\nabla \mathcal{U}(Q^\dag) \otimes \nabla w\right)\right) \\
&= \int_\Omega \left( \nabla \mathcal{U}(Q^\dag) \otimes \nabla w - P_{\mathcal{K}}\left(\nabla \mathcal{U}(Q^\dag) \otimes \nabla w\right) \right)\cdot \left(Q_h - P_{\mathcal{K}}\left(\nabla \mathcal{U}(Q^\dag) \otimes \nabla w\right)\right)\\
&~\quad - \int_\Omega \nabla \mathcal{U}(Q^\dag) \otimes \nabla w \cdot \left(Q_h - P_{\mathcal{K}}\left(\nabla \mathcal{U}(Q^\dag) \otimes \nabla w\right)\right)\\
&\le - \int_\Omega \nabla \mathcal{U}(Q^\dag) \otimes \nabla w \cdot \left(Q_h - P_{\mathcal{K}}\left(\nabla \mathcal{U}(Q^\dag) \otimes \nabla w\right)\right)\\
&\le \left| \int_\Omega \nabla \mathcal{U}(Q^\dag) \otimes \nabla w \cdot (Q^\dag - Q_h)\right|\\
&= \left| \int_\Omega (Q^\dag - Q_h) \nabla \mathcal{U}(Q^\dag) \cdot \nabla w\right|,
\end{align*}
which finishes the proof.
\end{proof}

To prove Theorem \ref{nu21***} we need the following auxiliary result.

\begin{lemma}\label{auxi2}
The estimate
$$\mathcal{J}^\delta_h (Q^\dag) \le C\left (\delta^2 + \sigma_h^2(Q^\dag) + \gamma^2_h \big( \mathcal{U}(Q^\dag) \big)\right)$$
holds.
\end{lemma}

\begin{proof}
The stated inequality follows from an argument which has included in the proof in Theorem \ref{convergence1} and therefore omitted.
\end{proof}

\begin{proof}[Proof of Theorem \ref{nu21***}]
Since $Q_h$ is the solution of the problem $\big(\mathcal{P}^{\rho,\delta}_h \big)$, we have that
\begin{align}
\mathcal{J}^\delta_h \big( Q_h \big) + \rho \big \|
Q_h  \big\|^2_{{L^2(\Omega)}^{d\times d}} &\le \mathcal{J}^\delta_h
\big(Q^\dag\big) + \rho \big\| Q^\dag  \big\|^2_{{L^2(\Omega)}^{d\times d}} \nonumber \\
&\le C\left (\delta^2 + \sigma_h^2(Q^\dag) + \gamma^2_h \big( \mathcal{U}(Q^\dag) \big)\right) + \rho \big\| Q^\dag  \big\|^2_{{L^2(\Omega)}^{d\times d}}, \label{27-4-15:ct1}
\end{align}
by Lemma \ref{auxi2}. Thus, we get
\begin{align}
\mathcal{J}_h^\delta \big( Q_h \big) &+ \rho \big\|
Q_h - Q^\dag  \big\|^2_{{L^2(\Omega)}^{d\times d}}
\nonumber\\
&\le C\left(\delta^2 + \sigma_h^2(Q^\dag) + \gamma_h^2 \big( \mathcal{U}(Q^\dag) \big) \right) \notag \\
&~\quad + \rho\left(\big\| Q^\dag
\big\|^2_{{L^2(\Omega)}^{d\times d}} - \big\| Q_h
\big\|^2_{{L^2(\Omega)}^{d\times d}}
+ \big\| Q_h - Q^\dag \big\|^2_{{L^2(\Omega)}^{d\times d}} \right)\nonumber\\
&= C\left (\delta^2 + \sigma_h^2(Q^\dag) + \gamma^2_h \big( \mathcal{U}(Q^\dag) \big)\right) + 2\rho \big\< Q^\dag  , Q^\dag -
Q_h \big\>_{{L^2(\Omega)}^{d \times d}}. \label{mqq11}
\end{align}
By \eqref{inf1}, \eqref{moi14***} and \eqref{29-6-15ct1}, we have with $w$ from \eqref{25-9-15ct1}
\begin{align}\label{28-9-15ct1}
\big\< Q^\dag  , Q^\dag - Q_h \big\>_{{L^2(\Omega)}^{d \times d}} &= \big\< Q^\dag  , Q^\dag - Q_h \big\>_{\big({{L^{\infty}(\Omega)}^{d \times d}}^*, {L^{\infty}(\Omega)}^{d \times d}\big)} \notag \\
&= \big \langle {\mathcal{U}'(Q^\dag)}^* w^*, Q^\dag - Q_h \big \rangle_{\big({{L^{\infty}(\Omega)}^{d \times d}}^*, {L^{\infty}(\Omega)}^{d \times d}\big)} \notag\\
&= \langle w^*, \mathcal{U}'(Q^\dag) (Q^\dag - Q_h) \rangle_{\big({H^1_0(\Omega)}^*, H^1_0(\Omega)\big)} \notag\\
&= \int_\Omega Q^\dag \nabla \mathcal{U}'(Q^\dag) (Q^\dag - Q_h) \cdot  \nabla w \notag\\
&= - \int_\Omega (Q^\dag - Q_h) \nabla \mathcal{U}(Q^\dag) \cdot  \nabla w,
\end{align}
here we used the equation \eqref{ct10}. Hence by (\ref{4/6:m4}), we get
\begin{align*}
\big\< Q^\dag  , Q^\dag - Q_h \big\>_{{L^2(\Omega)}^{d \times d}} &= \int_\Omega  Q_h
\nabla \mathcal{U} ( Q^\dag ) \cdot \nabla  w
- \int_\Omega f w \\
&= \int_\Omega  Q_h  \nabla \mathcal{U} ( Q^\dag ) \cdot
\nabla  w - \int_\Omega  Q_h \nabla \mathcal{U}
( Q_h ) \cdot \nabla  w \\
&= \int_\Omega  Q_h \nabla \big(\mathcal{U} ( Q^\dag ) -
\mathcal{U} ( Q_h )\big) \cdot \nabla  w \\
&= \int_\Omega  Q_h \nabla \big(\mathcal{U} (Q^\dag)-\Pi_hz^{\delta} \big) \cdot \nabla  w \\
&~\quad + \int_\Omega
Q_h \nabla \big(\mathcal{U}_h ( Q_h) - \mathcal{U}(Q_h)\big) \cdot \nabla  w\\
&~\quad + \int_\Omega
Q_h \nabla \big(\Pi_hz^{\delta} - \mathcal{U}_h (Q_h)\big) \cdot \nabla  w \\
&:= S_1 + S_2 + S_3.
\end{align*}
We have that
\begin{align}\label{24-6-15ct1}
\big\| \mathcal{U} (Q^\dag) -  \Pi_hz^{\delta} \big\|_{H^1 (\Omega)} &\le
\big\|  \Pi_h \big (\mathcal{U} (Q^\dag) - z^{\delta} \big) \big\|_{H^1 (\Omega)} +
\big\| \mathcal{U} (Q^\dag) - \Pi_h \mathcal{U} (Q^\dag) \big\|_{H^1 (\Omega)} \notag\\
&\le \big\| \Pi_h \big\|_{\mathcal {L} (H^1 (\Omega),
H^1(\Omega))} \big\|\mathcal{U} (Q^\dag) - z^{\delta} \big\|_{H^1
(\Omega)} + \big\| \mathcal{U} (Q^\dag) - \Pi_h \mathcal{U} (Q^\dag) \big\|_{H^1 (\Omega)} \notag\\
& \le \max \big(1, \| \Pi_h \|_{\mathcal{L}(H^1(\Omega), H^1(\Omega))} \big) \left(\delta + \gamma_h \big( \mathcal{U}(Q^\dag) \big)\right) \notag\\
&= C\left(\delta + \gamma_h \big( \mathcal{U}(Q^\dag) \big)\right).
\end{align}
Thus we obtain
\begin{align*}
S_1 &:= \int_\Omega  Q_h \nabla \big(\mathcal{U} (Q^\dag)-\Pi_hz^{\delta} \big) \cdot \nabla  w \\
&\le C \big \| \mathcal{U} (Q^\dag)-\Pi_hz^{\delta} \big \|_{H^1(\Omega)} \\
&\le C\left(\delta + \gamma_h \big( \mathcal{U}(Q^\dag) \big)\right).
\end{align*}
We deduce from \eqref{4/6:m4} and \eqref{10/4:ct1} that
$$\int_\Omega
Q_h \nabla \big(\mathcal{U}_h ( Q_h) - \mathcal{U}(Q_h)\big) \cdot \nabla v_h =0$$
for all $v_h \in \mathcal{V}^1_h.$ Therefore, we obtain
\begin{align*}
S_2 &:= \int_\Omega Q_h \nabla \big(\mathcal{U}_h ( Q_h) - \mathcal{U}(Q_h)\big) \cdot \nabla  w \\
&= \int_\Omega Q_h \nabla \big(\mathcal{U}_h ( Q_h) - \mathcal{U}(Q_h)\big) \cdot \nabla ( w - \Pi_h  w) \\
&\le \big\| Q_h \nabla \big(\mathcal{U}_h ( Q_h) - \mathcal{U}(Q_h)\big) \big\|_{L^2(\Omega)} \big\| \nabla ( w - \Pi_h  w) \big\|_{L^2(\Omega)} \\
&\le C \big\|  w - \Pi_h  w\big\|_{H^1(\Omega)} \\
&\le C \gamma_h(w).
\end{align*}
Since $0 \preceq Q_h (x) \in \mathcal{S}_d$ for a.e. $x \in \Omega$, the root ${Q_h (x)}^{1/2}$ is well defined. Furthermore,
\begin{align*}
Q_h(x) \nabla \big(\Pi_hz^{\delta}(x) - \mathcal{U}_h(Q_h)(x)\big) &\cdot \nabla  w(x) \\
&= {Q_h}^{1/2}(x) \nabla \big(\Pi_hz^{\delta}(x) -\mathcal{U}_h(Q_h)(x)\big)
\cdot {Q_h}^{1/2}(x) \nabla  w(x).
\end{align*}
Then applying the Cauchy-Schwarz inequality, we have
\begin{align*}
S_3 &:= \int_\Omega
Q_h \nabla \big(\Pi_hz^{\delta} - \mathcal{U}_h (Q_h)\big) \cdot \nabla  w \\
&\le \left(\int_\Omega Q_h \nabla \big(
\mathcal{U}_h(Q_h) -\Pi_hz^{\delta} \big)\cdot \nabla \big( \mathcal{U}_h(Q_h) -\Pi_hz^{\delta} \big)\right)^{1/2} \left(\int_\Omega
Q_h \nabla  w \cdot \nabla  w \right)^{1/2}\\
&\le C \left(\int_\Omega  Q_h \nabla
\big(\mathcal{U}_h ( Q_h )- \Pi_hz^{\delta} \big) \cdot \nabla
\big(\mathcal{U}_h ( Q_h )- \Pi_hz^{\delta} \big) \right)^{1/2}.
\end{align*}
Using Young's inequality, we obtain
\begin{align*}
S_3 \le C \rho
+\frac{1}{4\rho}\int_\Omega Q_h \nabla
\big(\mathcal{U}_h ( Q_h )- \Pi_hz^{\delta} \big)
\cdot \nabla \big(\mathcal{U}_h (Q_h )- \Pi_hz^{\delta} \big).
\end{align*}
Therefore, we arrive at
\begin{align*}
2\rho \big\< Q^\dag  , Q^\dag - Q_h \big\>_{{L^2(\Omega)}^{d \times d}} \le C\left(\delta^2 + \gamma_h^2 \big( \mathcal{U}(Q^\dag) \big) + \gamma^2_h (w) + \rho^2 \right) + \frac{1}{2}\mathcal{J}_h^\delta \big( Q_h \big).
\end{align*}
Combining this with \eqref{mqq11}, we infer that
\begin{align*}
\frac{1}{2}\mathcal{J}_h^\delta \big( Q_h \big) &+ \rho \big\|
Q_h - Q^\dag  \big\|^2_{{L^2(\Omega)}^{d\times d}} \le C\left(\delta^2 + \sigma_h^2(Q^\dag) + \gamma_h^2 \big( \mathcal{U}(Q^\dag) \big) + \gamma^2_h (w) + \rho^2 \right).
\end{align*}
Now we have
\begin{align*}
\dfrac{\kappa}{4} \big\|\mathcal{U}_h( Q_h ) - \mathcal{U}(Q^\dag)
\big\|^2_{H^1(\Omega)}
&\le \dfrac{\kappa}{2} \big\|\mathcal{U}_h( Q_h ) - \Pi_hz^{\delta} \big\|^2_{H^1(\Omega)} + \dfrac{\kappa}{2} \big\| \Pi_hz^{\delta} - \mathcal{U}(Q^\dag)\big\|^2_{H^1(\Omega)}\\
&\le \frac{1}{2}\mathcal{J}_h^\delta \big( Q_h \big) +\kappa C^2\left(\delta^2 + \gamma^2_h \big( \mathcal{U}(Q^\dag) \big)\right),
\end{align*}
by (\ref{coercivity}) and \eqref{24-6-15ct1}. Thus, we arrive at \eqref{29-6-15ct3},
which finishes the proof.
\end{proof}

\section{Gradient-projection algorithm} \label{iterative}

For the numerical solution we here use the gradient-projection algorithm of \cite{Enyi}. We note that many other efficient solution methods are available, see for example \cite{KeungZoz2000}.

We consider the finite dimensional space $\mathcal{V}_h$
defined by (\ref{18/10/14:ct1}). Let $\underline{C}_h$ and $\overline{C}_h$
be positive constants such that
\begin{align}\label{18/10/14:ct2}
\underline{C}_h \| H \|_{{L^2(\Omega)}^{d \times d}} \le
\| H \|_{{L^\infty(\Omega)}^{d \times d}} \le
\overline{C}_h \| H \|_{{L^2(\Omega)}^{d \times d}}
\end{align}
for all $H \in \mathcal{V}_h$.

The following results are useful.

\begin{lemma} \label{Lip*}
The discrete coefficient-to-solution operator $
\mathcal{U}_h $ 
is Lipschitz continuous on $\mathcal{Q}_{ad} \cap \mathcal{V}_h$ in the ${L^2(\Omega)}^{d \times d}$-norm
with a Lipschitz constant
\begin{align*}
\frac{\overline{C}_h d}{\kappa^2}
\big\|f\big\|_{L^2(\Omega)}.
\end{align*}
\end{lemma}
\begin{proof}
For all $M, N \in \mathcal{Q}_{ad} \cap \mathcal{V}_h$
it follows from (\ref{10/4:ct1}) that
\begin{align*}
\int_\Omega M\nabla \mathcal{U}_h(M) \cdot \nabla v_h &=
\int_\Omega fv_h \\
&= \int_\Omega N\nabla \mathcal{U}_h(N) \cdot \nabla v_h
\end{align*}
for all $v_h\in \mathcal{V}^1_h$. Thus
\begin{align*}
\int_\Omega M\nabla (\mathcal{U}_h(M) - \mathcal{U}_h(N))
\cdot \nabla v_h
= \int_\Omega (N-M)\nabla \mathcal{U}_h(N) \cdot
\nabla v_h.
\end{align*}
Choosing $v_h = \mathcal{U}_h(M) - \mathcal{U}_h(N)$,
by (\ref{coercivity}), we have
\begin{align*}
\kappa \big\| \mathcal{U}_h(M) - \mathcal{U}_h(N) \big
\|^2_{H^1(\Omega)} & \le d\| M - N \|_{{L^\infty(\Omega)}^{d \times d}}
\| \mathcal{U}_h(N) \|_{H^1(\Omega)}
\| \mathcal{U}_h(M) - \mathcal{U}_h(N) \|_{H^1(\Omega)}.
\end{align*}
Therefore, from (\ref{18/5:ct1}) and (\ref{18/10/14:ct2}) we arrive at
\begin{align} \label{18/10/14:ct3}
\| \mathcal{U}_h(M) - \mathcal{U}_h(N) \|_{H^1(\Omega)}
& \le \frac{\overline{C}_h d}{\kappa^2}
\big\|f\big\|_{L^2(\Omega)}
\| M - N \|_{{L^2(\Omega)}^{d \times d}}.
\end{align}
This finishes the proof.
\end{proof}

\begin{lemma} \label{Lip}
The objective functional $ \Upsilon^{\rho,\delta}_h $
of $\big(\mathcal{P}^{\rho,\delta}_h \big)$ has the
property that the gradient is Lipschitz continuous on
$\mathcal{Q}_{ad} \cap \mathcal{V}_h$ in the ${L^2(\Omega)}^{d \times d}$-norm with a Lipschitz constant
\begin{align*}
L_h := 2 \overline{C}_h d
\left( \frac{\overline{C}_h d}{\kappa^3}
\big\|f\big\|_{L^2(\Omega)}^2 + \rho |\Omega|^{1/2} \right).
\end{align*}
In other word, the estimate
$$\| {\Upsilon^{\rho,\delta}_h}' (M) -
{\Upsilon^{\rho,\delta}_h}' (N)
\|_{\mathcal{L}({L^2(\Omega)}^{d \times d}, {R})}
\le L_h \| M - N \|_{{L^2(\Omega)}^{d \times d}}$$
is satisfied for all
$M, N \in \mathcal{Q}_{ad} \cap \mathcal{V}_h$.
\end{lemma}

\begin{proof}
Since any norm on $\mathcal{V}_h$ is equivalent,
$\mathcal{U}_h$ is Fr\'echet differentiable on the set
$\mathcal{Q}_{ad} \cap \mathcal{V}_h$ in the
${L^\infty(\Omega)}^{d\times d}$-norm and thus in the
${L^2(\Omega)}^{d\times d}$-norm. For all
$M, N \in \mathcal{Q}_{ad} \cap \mathcal{V}_h$ and
$H \in \mathcal{V}_h$, in view of (\ref{17/10/14:ct6}), we get
\begin{align*}
\left| {\Upsilon^{\rho,\delta}_h}' (M)H
- {\Upsilon^{\rho,\delta}_h}' (N)H \right|
&= \Big| \int_{\Omega} H \nabla \mathcal{U}_h(N)
\cdot \nabla \mathcal{U}_h(N) -
\int_{\Omega} H \nabla \mathcal{U}_h(M)
\cdot \nabla \mathcal{U}_h(M)\\
&~\quad + 2\rho \int_\Omega H \cdot M
- 2\rho \int_\Omega H \cdot N\Big|.
\end{align*}
Thus
\begin{align*}
\Big| {\Upsilon^{\rho,\delta}_h}' (M)H
&- {\Upsilon^{\rho,\delta}_h}' (N)H \Big| \\
&= \Big| \int_{\Omega} H \nabla (\mathcal{U}_h(N) - \mathcal{U}_h(M))
\cdot \nabla(\mathcal{U}_h(N) + \mathcal{U}_h(M))
+ 2\rho \int_\Omega H \cdot (M - N) \Big|\\
&\le \left(\int_{\Omega} |H \nabla (\mathcal{U}_h(N) - \mathcal{U}_h(M))|^2\right)^{1/2}
\left(\int_{\Omega} | \nabla (\mathcal{U}_h(N)
+ \mathcal{U}_h(M))|^2\right)^{1/2}\\
&~\quad + 2\rho \left(\int_{\Omega} H \cdot H\right)^{1/2}
\left(\int_{\Omega} (M - N) \cdot (M - N)\right)^{1/2}\\
&\le d \| H \|_{{L^\infty(\Omega)}^{d \times d}}
\| \mathcal{U}_h(N) - \mathcal{U}_h(M) \|_{H^1(\Omega)}
\| \mathcal{U}_h(N) + \mathcal{U}_h(M) \|_{H^1(\Omega)}\\
&~\quad + 2\rho d \| H \|_{{L^\infty(\Omega)}^{d \times d}}  |\Omega|^{1/2}
\| M - N \|_{{L^2(\Omega)}^{d \times d}}.
\end{align*}
From the estimates (\ref{18/5:ct1}), (\ref{18/10/14:ct2}) and (\ref{18/10/14:ct3}) we now get
\begin{align*}
\Big| {\Upsilon^{\rho,\delta}_h}' (M)H &-
{\Upsilon^{\rho,\delta}_h}' (N)H \Big| \\
&\le d \| H \|_{{L^\infty(\Omega)}^{d \times d}}
\| \mathcal{U}_h(N) - \mathcal{U}_h(M) \|_{H^1(\Omega)}
\left(\| \mathcal{U}_h(N) \|_{H^1(\Omega)} +
\| \mathcal{U}_h(M) \|_{H^1(\Omega)}\right)\\
&~\quad + 2\rho d \| H \|_{{L^\infty(\Omega)}^{d \times d}}  |\Omega|^{1/2}
\| M - N \|_{{L^2(\Omega)}^{d \times d}}\\
&\le 2 \overline{C}_h d
\left( \frac{\overline{C}_h d}{\kappa^3}
\big\|f\big\|_{L^2(\Omega)}^2 + \rho |\Omega|^{1/2} \right)
\| H \|_{{L^2(\Omega)}^{d \times d}} \| M - N \|_{{L^2(\Omega)}^{d \times d}}.
\end{align*}
The lemma is proved.
\end{proof}

\begin{lemma} [\cite{Enyi}]\label{cite}
Let $X$ be a non-empty, closed and convex subset of a
Hilbert space $\mathcal{X}$ and
$\mathfrak{F}: X \to {R}$ be a convex Fr\'echet differentiable functional
with the gradient $\nabla \mathfrak{F}$ being $L$-Lipschitzian.
Assume that the problem
\begin{align} \label{21/10/2014:ct1}
\min_{x \in X} \mathfrak{F}(x)
\end{align}
is consistent and let $S$ denote its solution set.
Let ${(\alpha_m)}_m, {(\beta_m)}_m$ and ${(\gamma_m)}_m$
be real sequences satisfying ${(\alpha_m)}_m \subset (0,1),
~ \overline{{(\beta_m)}_m} \subset (0,1),
~ {(\gamma_m)}_m \subset (0, L/2)$ and the following additional condition
$$\lim_{m} \alpha_m = 0,~ \sum_{m = 1}^{\infty} \alpha_m
= \infty ~ \mbox{and} ~ 0< \liminf_{m} \gamma_m
\le \limsup_m \gamma_m < L/2.$$
Then, for any given $x^* \in X$ the iterative sequence
${(x_m)}_m$ is generated by $x_1 \in X$,
\begin{align}\label{25-9-15ct2}
x_{m+1} := (1 - \beta_m) x_m + \beta_m P_X
\left(x_m - \gamma_m \nabla \mathfrak{F} (x_m)\right)
+ \alpha_m (x^* - x_m)
\end{align}
converges strongly to the minimizer $x^\dag = P_S x^*$
of the problem (\ref{21/10/2014:ct1}).
\end{lemma}

To identify a stopping criterion for the iteration \eqref{25-9-15ct2} we adopt the following result.

\begin{lemma}\label{stopping criterion}
Let $X$ be a non-empty, closed and convex subset of a
Hilbert space $\mathcal{X}$ and
$\mathfrak{F}: X \to {R}$ be a convex Fr\'echet differentiable functional
with the gradient $\nabla \mathfrak{F}$. Assume that the problem
\begin{align} \label{25-9-15ct3}
\min_{x \in X} \mathfrak{F}(x)
\end{align}
is consistent. Then $x^\dag$ is a solution to \eqref{25-9-15ct3} if and only if the equation
$$x^\dag = P_X \left( x^\dag - \gamma
\nabla \mathfrak{F} (x^\dag)\right) $$
holds, where $\gamma$ is an arbitrary positive constant.
\end{lemma}

\begin{proof}
Since $\mathfrak{F}$ is convex differentiable, we have for all $\gamma>0$ that
\begin{align*}
x^\dag \mbox{solves \eqref{25-9-15ct3}} &\Leftrightarrow \left\langle \gamma\nabla \mathfrak{F} (x^\dag), x-x^\dag \right\rangle _{\mathcal{X}} \ge 0 ~\quad \mbox{for all~} x \in X\\
&\Leftrightarrow \left\langle  \left( x^\dag - \gamma\nabla \mathfrak{F} (x^\dag) \right)  - x^\dag, x-x^\dag\right\rangle _{\mathcal{X}} \le 0 ~\quad \mbox{for all~} x \in X \\
& \Leftrightarrow x^\dag = P_X \left( x^\dag - \gamma
\nabla \mathfrak{F} (x^\dag)\right),
\end{align*}
which finishes the proof.
\end{proof}

Now we state the main result of this section on the strong
convergence of iterative solutions to that of our identification problem.

\begin{theorem}\label{algorithm}
Let $\left(\mathcal{T}_{h_n}\right)_n$ be a sequence of triangulations with
mesh sizes $\left(h_n\right)_n$. For any positive sequence
$\left(\delta_n\right)_n$, let $\rho_n := \rho\left(\delta_n\right)$
be such that
$$\rho_n\rightarrow 0, ~\frac{\delta_n^2}{\rho_n}
\rightarrow 0, ~\frac{\sigma^2_{h_n}(Q^\dag)}{\rho_n}
\rightarrow 0 \mbox{~and~} \frac{\gamma^2_{h_n}(Q^\dag)}{\rho_n}
\rightarrow 0$$
and $\left(z^{\delta_n}\right)_n$ be observations satisfying
$\big\| \overline{u} - z^{\delta_n}
\big\|_{H^1_0(\Omega)} \le \delta_n$.

Moreover, for any fixed $n$ let
${(\alpha^n_m)}_m, {(\beta^n_m)}_m$ and ${(\gamma^n_m)}_m$
be real sequences satisfying
\begin{quote}
${(\alpha^n_m)}_m \subset (0,1), ~ \overline{{(\beta^n_m)}_m} \subset (0,1),
~{(\gamma^n_m)}_m \subset (0, L_{h_n}/2)$,

$\lim_{m} \alpha^n_m = 0,~ \sum_{m = 1}^{\infty} \alpha^n_m = \infty$ and

$0< \liminf_m \gamma^n_m \le \limsup_m \gamma^n_m < L_{h_n}/2$
\end{quote}
with
\begin{align}\label{28-4-15-ct1}
L_{h_n}:= 2 \overline{C}_{h_n} d \left( \frac{\overline{C}_{h_n} d}{\kappa^3} \big\|f\big\|_{L^2(\Omega)}^2 + \rho_n  |\Omega|^{1/2}\right).
\end{align}
Let $Q^*$ be a prior estimate of the sought matrix $Q^\dag$ and let ${(Q^n_m)}_m$ be the sequence of iterates generated by
\begin{align}\label{30-4-15:ct1}
\begin{split}
&Q^n_0 \in \mathcal{Q}_{ad} \cap \mathcal{V}_{h_n}\\
&Q^n_{m} := (1 - \beta^n_{m-1}) Q^n_{m-1} + \alpha^n_{m-1} (Q^* - Q^n_{m-1}) \\
&\phantom{xxx} + \beta^n_{m-1} P_{\mathcal{Q}_{ad} \cap \mathcal{V}_{h_n}}
\big(Q^n_{m-1} - \gamma^n_{m-1} \big( \nabla \Pi_{h_n}
z^{\delta_n} \otimes \nabla \Pi_{h_n }z^{\delta_n}
-\nabla \mathcal{U}_{h_n}(Q^n_{m-1}) \otimes
\nabla \mathcal{U}_{h_n}(Q^n_{m-1}) \\
&\phantom{xxx} + 2\rho_n Q^n_{m-1}\big)\big).
\end{split}
\end{align}
Then ${(Q^n_m)}_m$ converges strongly to the unique
minimizer $ Q^{\rho_n,\delta_n}_{h_n}$ of
$\left( \mathcal{P}^{\rho_n,\delta_n}_{h_n} \right)$,
$$\lim_m \| Q^n_m - Q^{\rho_n,\delta_n}_{h_n}
\|_{{L^2(\Omega)}^{d \times d}} = 0. $$
Furthermore, ${(Q^n_m)}_m^n$ converges strongly to the minimum norm solution $Q^\dag$ of the
identification problem,
$$\lim_n \big( \lim_m \| Q^n_m - Q^\dag
\|_{{L^2(\Omega)}^{d \times d}} \big) = 0. $$
\end{theorem}

\begin{proof}
Since
$$\nabla \Upsilon^{\rho,\delta}_h (Q) = \nabla
 \Pi_hz^{\delta} \otimes \nabla  \Pi_hz^{\delta} -
\nabla \mathcal{U}_h(Q) \otimes
\nabla \mathcal{U}_h(Q) + 2\rho Q$$
for all $Q \in \mathcal{Q}_{ad} \cap \mathcal{V}_h$,
the conclusion of the theorem follows directly from
Theorem \ref{convergence1}, Lemma \ref{Lip} and Lemma \ref{cite}.
\end{proof}

\section{Numerical tests}\label{Numerical implement}

In this section we illustrate the theoretical result with two numerical examples. The first one is provided to compare with the numerical results obtained in \cite{Deckelnick}, while the second one aims to illustrate the discontinuous coefficient identification problem.

For this purpose we consider the Dirichlet problem
\begin{align}
-\divv (Q^\dag\nabla \overline{u}) &= f \mbox{ in } \Omega, \label{m1**}\\
\overline{u} &=0 \mbox{ on } {\partial \Omega} \label{qmict3**}
\end{align}
with $\Omega = \{ x = (x_1,x_2) \in {R}^2 ~|~ -1 < x_1, x_2 < 1\}$ and
\begin{align}\label{4-3-16ct1}
\overline{u}(x) &= (1-x^2_1)(1-x^2_2).
\end{align}
Now we divide the interval $(-1,1)$ into $\ell$ equal segments and so that the domain $\Omega = (-1,1)^2$ is divided into $2\ell^2$ triangles, where the diameter of each triangle is $h_{\ell} = \frac{\sqrt{8}}{\ell}$. In the minimization problem $\left(\mathcal{P}^{\rho,\delta}_h \right)$ we take $h=h_\ell$ and $\rho = \rho_\ell = 0.001h_\ell$. For observations with noise we assume that
$$z^{\delta_{\ell}} = \overline{u} + \dfrac{x_1}{\ell} + \dfrac{x_2}{\ell} \mbox{~and~} \Pi_{h_\ell} z^{\delta_{\ell}} = I^1_{h_{\ell}} \left( \overline{u} + \dfrac{x_1}{\ell} + \dfrac{x_2}{\ell}\right) $$
so that
$$\delta_\ell = \left\| z^{\delta_{\ell}} - \overline{u} \right\|_{H^1(\Omega)} = \left\| \dfrac{x_1}{\ell} + \dfrac{x_2}{\ell}\right\|_{H^1(\Omega)} = \sqrt{\frac{32}{3}} \frac{1}{\ell}= \dfrac{2}{\sqrt{3}} h_\ell.$$
The constants $\underline{q}$ and $\overline{q}$ in the definition of the set $\mathcal{K}$ are respectively chosen as 0.05 and 10. We use the gradient-projection algorithm which is described in Theorem \ref{algorithm} for computing the solution of the problem $\left(\mathcal{P}^{\rho,\delta}_h \right)$.

Note that in \eqref{28-4-15-ct1} $d=2$ and $\overline{C}_{h_{\ell}} = \frac{\ell}{\sqrt{2}}$, where for all $Q\in \mathcal{Q}_{ad}$ and $v \in H^1_0(\Omega)$
\begin{align*}
\|v\|^2_{H^1(\Omega)} &= \int_\Omega |\nabla v|^2 + \int_\Omega | v|^2\\
&\le \int_\Omega |\nabla v|^2 + \left( \sqrt{\frac{3}{2}} \right)^{(d+2)/2} |\Omega|^{1/d} \int_\Omega | \nabla v|^2\\
&\le \dfrac{1}{\underline{q}} \left( 1 + \left( \sqrt{\frac{3}{2}} \right)^{(d+2)/2} |\Omega|^{1/d}\right) \int_\Omega Q \nabla v \cdot \nabla v.
\end{align*}
So we can choose
\begin{align*}
\kappa = \dfrac{\underline{q}}{1 + \left( \sqrt{\frac{3}{2}} \right)^{(d+2)/2} |\Omega|^{1/d}}.
\end{align*}
As the initial matrix $Q_0$ in \eqref{30-4-15:ct1} we choose
$$Q_0 := \begin{pmatrix}
2&0\\
0&2
\end{pmatrix}.$$
The prior estimate is chosen with $Q^* := I^1_{h_\ell} Q^\dag$. Furthermore,
The sequences ${(\alpha_m)}_m, {(\beta_m)}_m$ and ${(\gamma_m)}_m$ are chosen with
$$\alpha_m = \frac{1}{100m}, ~~ \beta_m = \frac{100m\rho_{\ell}}{3m+1}, ~~\mbox{and}~ \gamma_m = \frac{100m\rho_{\ell}}{2m+1}.$$
Let $Q^{\rho_{\ell},\delta_{\ell}}_{h_{\ell}}$ denote the computed numerical matrix with respect to $\ell$ and the iteration \eqref{30-4-15:ct1}. According to the lemma \ref{stopping criterion}, the iteration was stopped if
$$\mbox{Tolerance~}:= \left\|Q^{\rho_{\ell},\delta_{\ell}}_{h_{\ell}} -  P_{\mathcal{Q}_{ad} \cap \mathcal{V}_{h_\ell}} \left(  Q^{\rho_{\ell},\delta_{\ell}}_{h_{\ell}} - \gamma_m \nabla \Upsilon^{\rho_{\ell},\delta_{\ell}}_{h_{\ell}} \big( Q^{\rho_{\ell},\delta_{\ell}}_{h_{\ell}} \big) \right)  \right\|_{{L^2(\Omega)}^{d\times d}} < 10^{-6}$$
with
$$\nabla \Upsilon^{\rho_{\ell},\delta_{\ell}}_{h_{\ell}} \big( Q^{\rho_{\ell},\delta_{\ell}}_{h_{\ell}} \big) = \nabla
 \Pi_{h_\ell}z^{\delta_\ell} \otimes \nabla \Pi_{h_\ell}z^{\delta_\ell} -
\nabla \mathcal{U}_{h_\ell} \big( Q^{\rho_{\ell},\delta_{\ell}}_{h_{\ell}} \big) \otimes
\nabla \mathcal{U}_{h_\ell} \big( Q^{\rho_{\ell},\delta_{\ell}}_{h_{\ell}} \big) + 2\rho_\ell Q^{\rho_{\ell},\delta_{\ell}}_{h_{\ell}}$$
or the number of iterations was reached 500.

\subsection{Example 1}\label{eg1}

We now assume that
\begin{align*}
Q^\dag (x) &= P_{\mathcal{K}}\left( \nabla \overline{u}(x) \otimes \nabla \overline{u}(x)\right).
\end{align*}

Let us denote $\eta(x) = 4(x_1^2(1-x_2^2)^2 + x_2^2(1-x_1^2)^2)$ and $P_{[\underline{q}, \overline{q}]} (\eta(x)) = \max\left(\underline{q}, \min(\eta(x), \overline{q})\right)$. A calculation shows
$$Q^\dag(x) = \begin{cases}
\underline{q}I_2 & \mbox{~if~} \eta(x)=0,\\
\underline{q}I_2 + \frac{P_{[\underline{q}, \overline{q}]} (\eta(x)) - \underline{q}}{\eta(x)} \nabla \overline{u}(x) \otimes \nabla \overline{u}(x) & \mbox{~if~} \eta(x) \neq 0.
\end{cases}$$
Then along with $\overline{u}$ given in the equation \eqref{4-3-16ct1} one can compute the right hand side $f$ in the equation  \eqref{m1**}.

The numerical results are summarized in Table \ref{b1}, where we present the refinement level $\ell$, regularization parameter $\rho_\ell$, mesh size $h_\ell$ of the triangulation, noise level $\delta_\ell$, number of iterations, value of tolerances, the final $L^2$ and $L^\infty$-error in the coefficients, the final $L^2$ and $H^1$-error in the states. Their experimental order of convergence (EOC) is presented in Table \ref{b11}, where
$$\mbox{EOC}_\Phi := \dfrac{\ln \Phi(h_1) - \ln \Phi(h_2)}{\ln h_1 - \ln h_2}$$
and $\Phi(h)$ is an error functional with respect to the mesh size $h$.

In Table \ref{b2} we present the numerical result for $\ell = 96$, where the value of tolerances, the final $L^2$ and $L^\infty$-error in the coefficients, the final $L^2$ and $H^1$-error in the states are displayed for each one hundred iteration. The convergence history given in Table \ref{b1}, Table \ref{b11} and Table \ref{b2} shows that the gradient-projection algorithm performs well for our identification problem.

\begin{figure}[H]
    \begin{center}
    \includegraphics[scale=0.27]{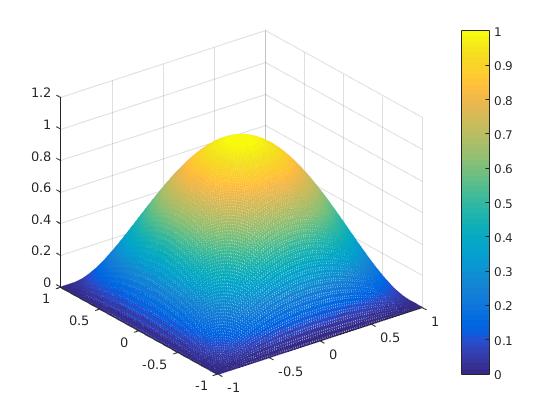}~~~~~~~~\includegraphics[scale=0.27]{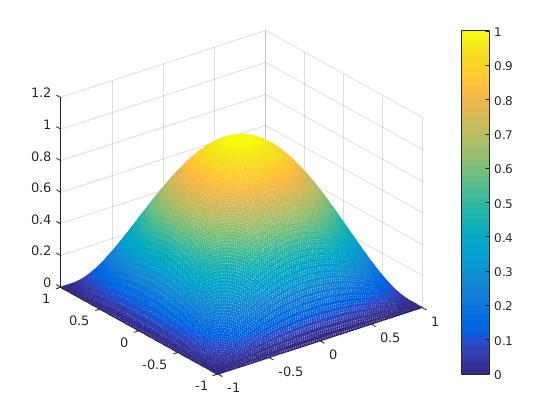} 			   ~~~~~~~~\includegraphics[scale=0.27]{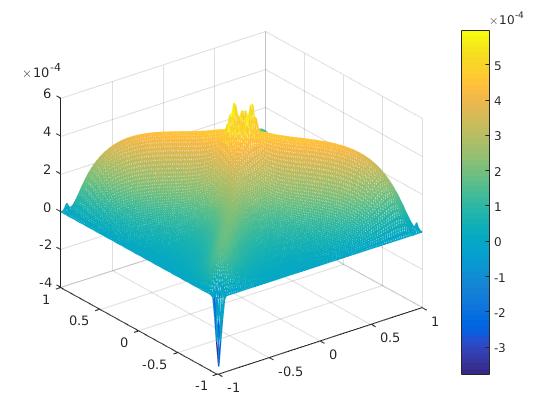}\\
    \end{center}
    \caption{Interpolation $I^1_{h_{\ell}} \overline{u}$, computed numerical state of the algorithm at the 500$^{\mbox{\tiny th}}$ iteration, and the difference to $I^1_{h_{\ell}}\overline{u}$.}
    \label{h1}
    \end{figure}

\begin{figure}[H]
    \begin{center}
    \includegraphics[scale=0.27]{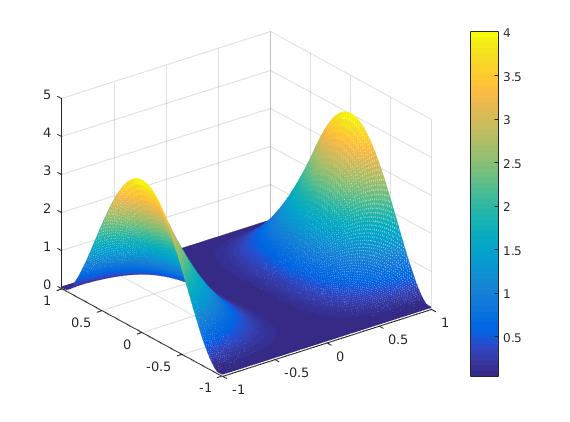}~~~~~~~~\includegraphics[scale=0.27]{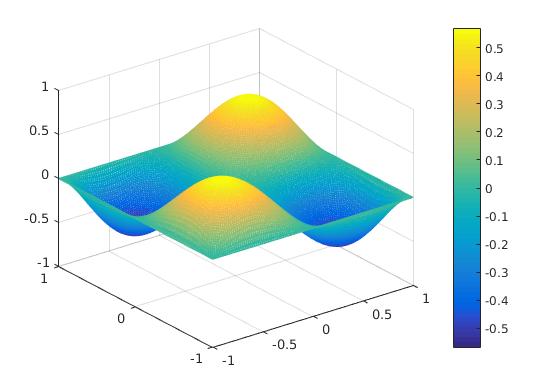} 			   ~~~~~~~~\includegraphics[scale=0.27]{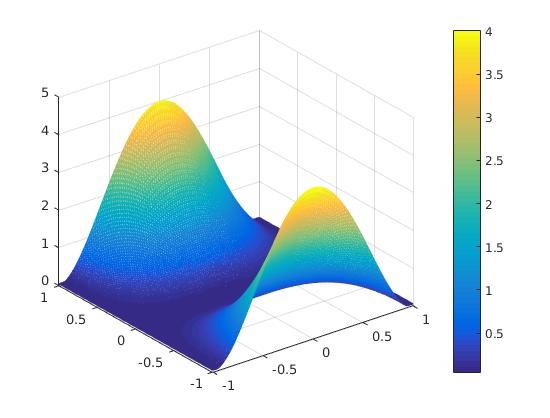}\\
    \end{center}
    \caption{Graphs of $I^1_{h_{\ell}} q^\dag_{11}, I^1_{h_{\ell}}q^\dag_{12}$ and $I^1_{h_{\ell}}q^\dag_{22}$.}
    \label{h3}
    \end{figure}

\begin{figure}[H]
    \begin{center}
    \includegraphics[scale=0.27]{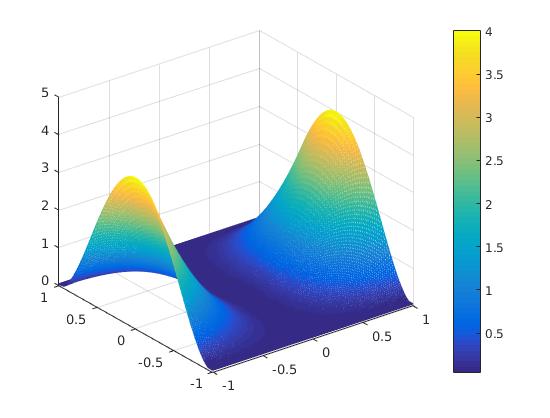}~~~~~~~~\includegraphics[scale=0.27]{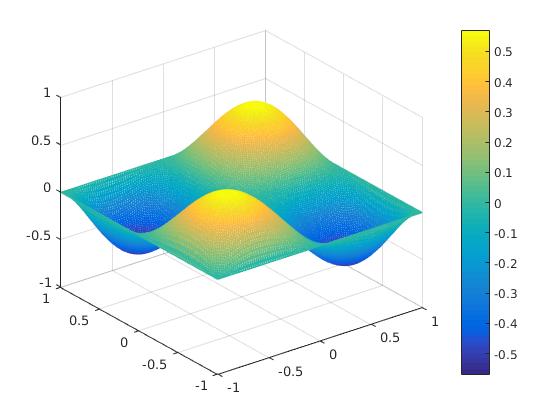} 			   ~~~~~~~~\includegraphics[scale=0.27]{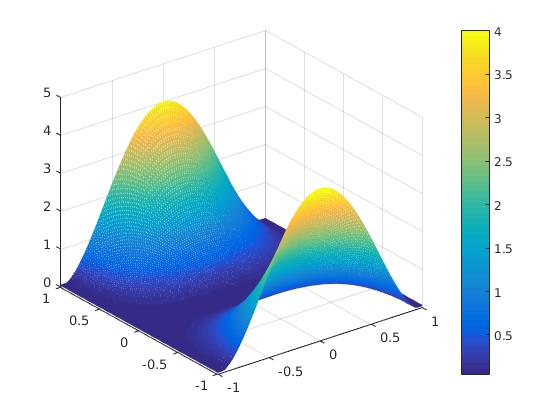}\\
    \end{center}
    \caption{Graphs of ${q^{\rho_{\ell},\delta_{\ell}}_{h_{\ell}}}_{11}, {q^{\rho_{\ell},\delta_{\ell}}_{h_{\ell}}}_{12}$ and ${q^{\rho_{\ell},\delta_{\ell}}_{h_{\ell}}}_{22}$.}
    \label{h4}
    \end{figure}

\begin{figure}[H]
    \begin{center}
    \includegraphics[scale=0.27]{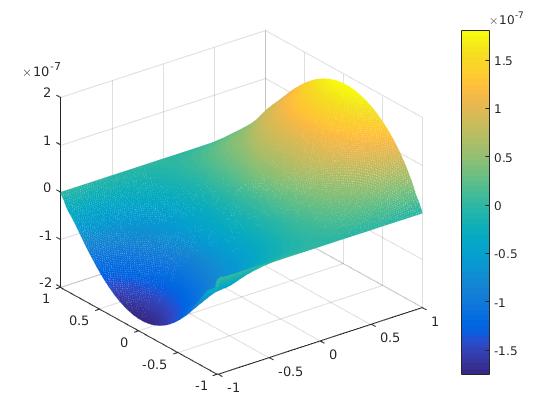}~~~~~~~~\includegraphics[scale=0.27]{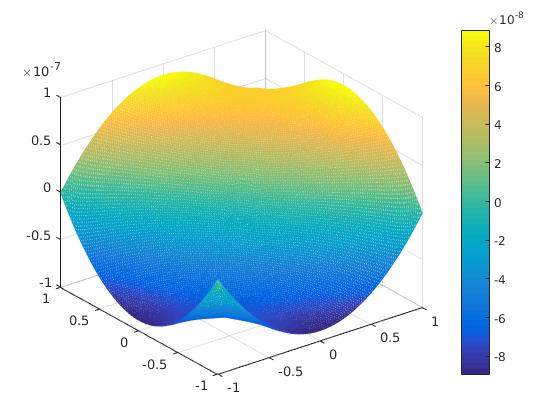} 			   ~~~~~~~~\includegraphics[scale=0.27]{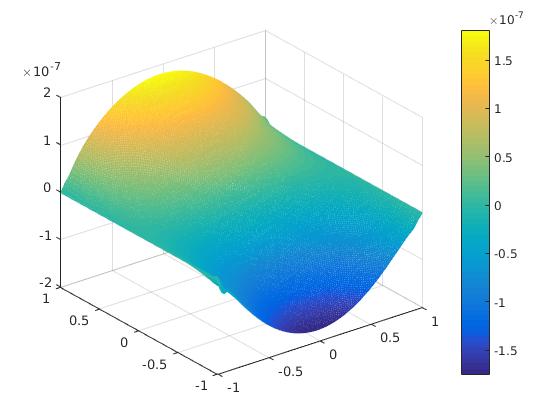}\\
    \end{center}
    \caption{Differences ${q^{\rho_{\ell},\delta_{\ell}}_{h_{\ell}}}_{11} - I^1_{h_{\ell}} q^\dag_{11}$, ${q^{\rho_{\ell},\delta_{\ell}}_{h_{\ell}}}_{12} - I^1_{h_{\ell}} q^\dag_{12}$ and ${q^{\rho_{\ell},\delta_{\ell}}_{h_{\ell}}}_{22} - I^1_{h_{\ell}} q^\dag_{22}$.}
    \label{h5}
    \end{figure}

All figures are presented here corresponding to $\ell = 96$. Figure \ref{h1} from left to right shows the graphs of the interpolation $I^1_{h_{\ell}} \overline{u}$, computed numerical state of the algorithm at the 500$^{\mbox{\tiny th}}$ iteration, and the difference to $I^1_{h_{\ell}}\overline{u}$.
We write
$$Q^\dag = \begin{pmatrix}
q^\dag_{11}&q^\dag_{12}\\
q^\dag_{12}&q^\dag_{22}
\end{pmatrix} ~~\mbox{and}~~ Q^{\rho_{\ell},\delta_{\ell}}_{h_{\ell}} = \begin{pmatrix}
{q^{\rho_{\ell},\delta_{\ell}}_{h_{\ell}}}_{11} & {q^{\rho_{\ell},\delta_{\ell}}_{h_{\ell}}}_{12}\\
{q^{\rho_{\ell},\delta_{\ell}}_{h_{\ell}}}_{12} & {q^{\rho_{\ell},\delta_{\ell}}_{h_{\ell}}}_{22}
\end{pmatrix}.$$
Figure \ref{h3} from left to right we display $I^1_{h_{\ell}} q^\dag_{11}, I^1_{h_{\ell}}q^\dag_{12}$ and $I^1_{h_{\ell}}q^\dag_{22}$. Figures \ref{h4} shows ${q^{\rho_{\ell},\delta_{\ell}}_{h_{\ell}}}_{11}, {q^{\rho_{\ell},\delta_{\ell}}_{h_{\ell}}}_{12}$ and ${q^{\rho_{\ell},\delta_{\ell}}_{h_{\ell}}}_{22}$. Figure \ref{h5} from left to right we display differences ${q^{\rho_{\ell},\delta_{\ell}}_{h_{\ell}}}_{11} - I^1_{h_{\ell}} q^\dag_{11}$, ${q^{\rho_{\ell},\delta_{\ell}}_{h_{\ell}}}_{12} - I^1_{h_{\ell}} q^\dag_{12}$ and ${q^{\rho_{\ell},\delta_{\ell}}_{h_{\ell}}}_{22} - I^1_{h_{\ell}} q^\dag_{22}$.

For the simplicity of the notation we denote by
\begin{align}
& \Gamma := \|Q^{\rho_{\ell},\delta_{\ell}}_{h_{\ell}} - I^1_{h_\ell} Q^\dag\|_{{L^2(\Omega)}^{d\times d}},~ \Delta := \|Q^{\rho_{\ell},\delta_{\ell}}_{h_{\ell}} - I^1_{h_\ell} Q^\dag\|_{{L^\infty(\Omega)}^{d\times d}}, \label{4-3-16ct2}\\
& \Sigma := \|\mathcal{U}_{h_{\ell}} \big( Q^{\rho_{\ell},\delta_{\ell}}_{h_{\ell}} \big) - I^1_{h_\ell}\overline{u}\|_{L^2(\Omega)} \mbox{~and~} \Xi := \|\mathcal{U}_{h_{\ell}} \big( Q^{\rho_{\ell},\delta_{\ell}}_{h_{\ell}} \big) - I^1_{h_\ell}\overline{u}\|_{H^1(\Omega)}.\label{4-3-16ct3}
\end{align}

\begin{table}[H]
\begin{center}
\begin{tabular}{|c|l|l|l|l|l|l|l|l|l|}
 \hline \multicolumn{10}{|c|}{ {\bf Convergence history} }\\
 \hline
$\ell$ &\scriptsize $\rho_\ell$ &\scriptsize $h_\ell$ &\scriptsize $\delta_\ell$ &\scriptsize {\bf Ite.} &\scriptsize {\bf Tol.} &\scriptsize  $\Gamma$ &\scriptsize $\Delta$ &\scriptsize $\Sigma$ &\scriptsize $\Xi$ \\
\hline
6   &4.7140e-4 &0.4714 &0.5443& 500& 0.0165 & 1.0481e-3 & 1.0381e-3  & 0.041001  &0.20961 \\
\hline
12  &2.3570e-4 &0.2357 &0.2722& 500& 0.0057 & 1.3471e-4  & 1.0825e-4  & 0.012848 &0.070352 \\
\hline
24  &1.1785e-4 &0.1179 &0.1361& 500& 0.0014 & 1.6826e-5  & 1.2084e-5  & 4.1855e-3  &0.023265 \\
\hline
48  &5.8926e-5 &0.0589 &0.0680& 500& 3.5653e-4 & 2.0962e-6  & 1.4577e-6  & 9.8725e-4  &6.3264e-3 \\
\hline
96  &2.9463e-5 &0.0295 &0.0340& 500& 8.9059e-5 & 2.6177e-7  & 1.8007e-7 & 2.6475e-4  &1.8936e-3 \\
\hline
\end{tabular}
\caption{Refinement level $\ell$, regularization parameter $\rho_\ell$, mesh size $h_\ell$ of the triangulation, noise level $\delta_\ell$, number of iterations, value of tolerances, errors  $\Gamma$, $\Delta$, $\Sigma$ and $\Xi$.}
\label{b1}
\end{center}
\end{table}

\begin{table}[H]
\begin{center}
\begin{tabular}{|c|c|c|c|c|}
 \hline \multicolumn{5}{|c|}{ {\bf Experimental order of convergence} }\\
 \hline
$\ell$ &\scriptsize {\bf EOC$_\Gamma$} &\scriptsize {\bf EOC$_\Delta$} &\scriptsize {\bf EOC$_\Sigma$} &\scriptsize {\bf EOC$_\Xi$}\\
\hline
6 & --& -- & -- &--\\
\hline
12 & 2.9598 & 3.2615 & 1.6741& 1.5750\\
\hline
24 & 3.0011 & 3.1632 & 1.6181 &1.5964 \\
\hline
48 & 3.0048 & 3.0513 & 2.0839 &1.8787\\
\hline
96 & 3.0014 & 3.0171 & 1.8988 &1.7403 \\
\hline
Mean of EOC & 2.9918 & 3.1233 & 1.8187 & 1.6976 \\
\hline
\end{tabular}
\caption{Experimental order of convergence between finest and coarsest
level for $\Gamma$, $\Delta$, $\Sigma$ and $\Xi$.}
\label{b11}
\end{center}
\end{table}

\begin{table}[H]
\begin{center}
\begin{tabular}{|c|l|l|l|l|l|}
 \hline \multicolumn{6}{|c|}{ {\bf Numerical result for $\ell = 96$}} \\
 \hline
{ \bf Iterations} &\scriptsize {\bf Tolerances} &\scriptsize  $\Gamma$ &\scriptsize $\Delta$ &\scriptsize $\Sigma$ &\scriptsize $\Xi$\\
\hline
100 & 0.1875 & 5.0541 & 4.5433 & 4.6756 & 27.2404\\
\hline
200 & 9.8522e-3 & 0.01224 & 0.58589 & 6.6362e-3 & 0.021169\\
\hline
300 & 8.9476e-5 & 5.8712e-7 & 6.0808e-7 & 2.6475e-4 & 1.8936e-3\\
\hline
400 & 8.9386e-5 & 3.0944e-7 & 2.4408e-7 & 2.6475e-4 & 1.8936e-3\\
\hline
500 & 8.9059e-5 & 2.6177e-7 & 1.8007e-7 & 2.6475e-4 & 1.8936e-3\\
\hline
 \end{tabular}
 \caption{Errors  $\Gamma$, $\Delta$, $\Sigma$ and $\Xi$ for $\ell = 96$.}
\label{b2}
\end{center}
\end{table}

\subsection{Example 2}\label{eg2}

We next assume that entries of the symmetric matrix $Q^\dag \in \mathcal{Q}_{ad}$ are discontinuous which are defined as
$$q^\dag_{11}(x) = \begin{cases}
3 & \mbox{~if~} x\in\Omega_{11}\\
1 & \mbox{~if~} x\in\Omega\setminus\Omega_{11}
\end{cases},~ q^\dag_{12}(x) = \begin{cases}
1 & \mbox{~if~} x\in\Omega_{12}\\
0 & \mbox{~if~} x\in\Omega\setminus\Omega_{12}
\end{cases} \mbox{~and~} q^\dag_{22}(x) = \begin{cases}
4 & \mbox{~if~} x\in\Omega_{22}\\
2 & \mbox{~if~} x\in\Omega\setminus\Omega_{22}
\end{cases}$$
with
\begin{align*}
&\Omega_{11} := \left\{ (x_1, x_2) \in \Omega ~\Big|~ |x_1| \le \frac{1}{2} \mbox{~and~} |x_2| \le \frac{1}{2} \right\},\\
&\Omega_{12} := \left\{ (x_1, x_2) \in \Omega ~\Big|~ |x_1| + |x_2| \le \frac{1}{2} \right\} \mbox{~and~}\\
&\Omega_{22} := \left\{ (x_1, x_2) \in \Omega ~\Big|~ x_1^2 + x_2^2 \le \left( \frac{1}{2}\right)^2 \right\}.
\end{align*}
Since the entries of the matrix $Q^\dag$ are discontinuous, the right hand side $f$ in the equation \eqref{m1**} is now given in the form of a load vector 
$$F = KU,$$
where $K=\left( k_{ij}\right)_{1\le i,j\le N_\ell}$ with
$$k_{ij} := \int_\Omega Q^\dag \nabla \phi_i \cdot \nabla \phi_j$$
and $\{\phi_1, \cdots, \phi_{N_\ell}\}$ being the basis for the approximating subspace $\mathcal{V}^1_{h_\ell}$, while the vector $U$ is the nodal values of the functional $\overline{u}$.

With the notation on errors $\Gamma$, $\Delta$, $\Sigma$ and $\Xi$ as in the equations \eqref{4-3-16ct2}-\eqref{4-3-16ct3} the numerical results of Example \ref{eg2} are summarized in Table \ref{b4}. For clarity we also present additionally the $H^1(\Omega)$-semi-norm error
$$\Lambda := \| \nabla \mathcal{U}_{h_{\ell}} \big( Q^{\rho_{\ell},\delta_{\ell}}_{h_{\ell}} \big) - \nabla I^1_{h_\ell}\overline{u}  \|_{L^2(\Omega)}.$$
For simplicity in Table \ref{b4} we do not restate the regularization parameter $\rho_\ell$, mesh size $h_\ell$ of the triangulation and noise level $\delta_\ell$ again, since they have been given in Table \ref{b1} of Example \ref{eg1}. 

The experimental order of convergence for $\Gamma$, $\Delta$, $\Sigma$, $\Xi$ and $\Lambda$ is presented in Table \ref{b5}. 

\begin{figure}[H]
\begin{center}
\includegraphics[scale=0.2]{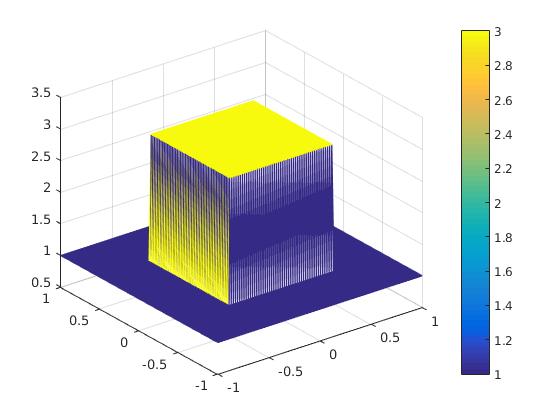}
\includegraphics[scale=0.2]{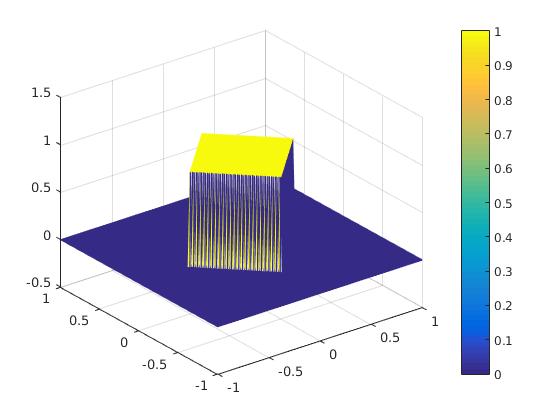}
\includegraphics[scale=0.2]{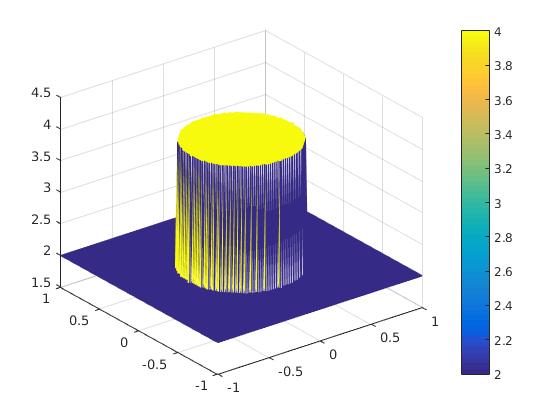}
\includegraphics[scale=0.2]{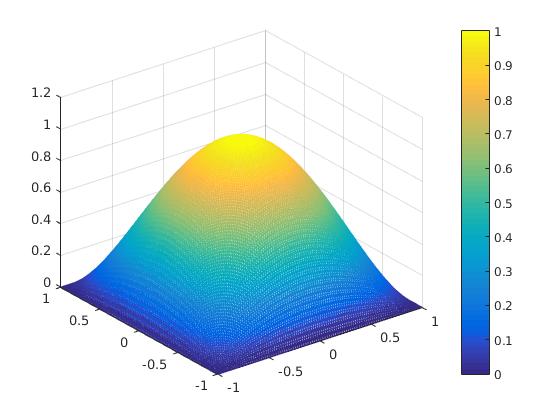}
\end{center}
\caption{Graphs of ${q^{\rho_{\ell},\delta_{\ell}}_{h_{\ell}}}_{11}, {q^{\rho_{\ell},\delta_{\ell}}_{h_{\ell}}}_{12}$, ${q^{\rho_{\ell},\delta_{\ell}}_{h_{\ell}}}_{22}$ and $\mathcal{U}_{h_{\ell}} \big( Q^{\rho_{\ell},\delta_{\ell}}_{h_{\ell}} \big)$ at the 500$^{\mbox{\tiny th}}$ iteration.}
\label{h6}
\end{figure}

\begin{figure}[H]
\begin{center}
\includegraphics[scale=0.2]{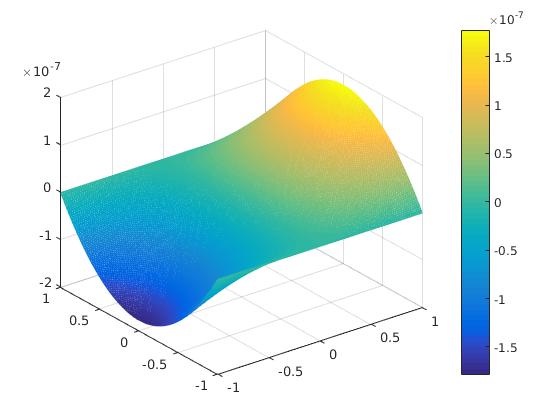}
\includegraphics[scale=0.2]{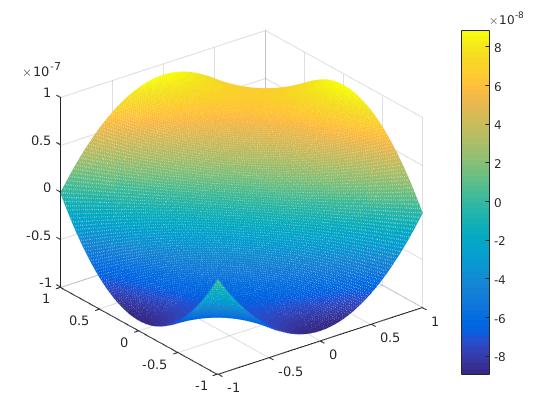}
\includegraphics[scale=0.2]{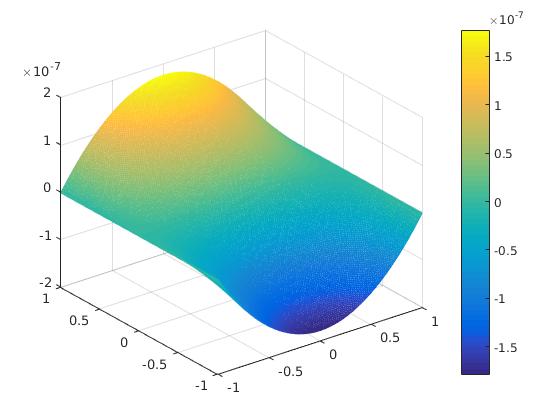}
\includegraphics[scale=0.2]{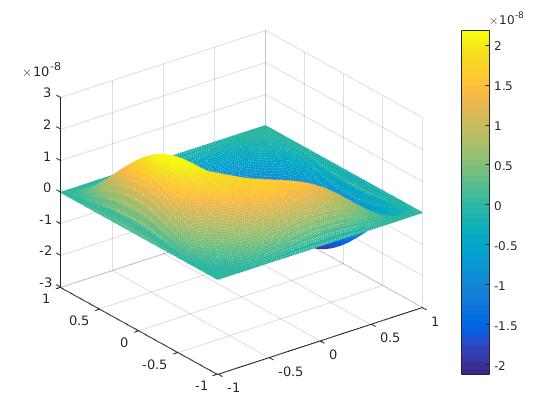}
\end{center}
\caption{Graphs of ${q^{\rho_{\ell},\delta_{\ell}}_{h_{\ell}}}_{11} - I^1_{h_{\ell}} q^\dag_{11}$, ${q^{\rho_{\ell},\delta_{\ell}}_{h_{\ell}}}_{12} - I^1_{h_{\ell}} q^\dag_{12}$, ${q^{\rho_{\ell},\delta_{\ell}}_{h_{\ell}}}_{22} - I^1_{h_{\ell}} q^\dag_{22}$ and $\mathcal{U}_{h_{\ell}} \big( Q^{\rho_{\ell},\delta_{\ell}}_{h_{\ell}} \big) - I^1_{h_\ell}\overline{u}$ at the 500$^{\mbox{\tiny th}}$ iteration.}
\label{h7}
\end{figure}

All figures are presented corresponding to $\ell = 96$. Figure \ref{h6} from left to right contains graphs of the entries  ${q^{\rho_{\ell},\delta_{\ell}}_{h_{\ell}}}_{11}, {q^{\rho_{\ell},\delta_{\ell}}_{h_{\ell}}}_{12}$, ${q^{\rho_{\ell},\delta_{\ell}}_{h_{\ell}}}_{22}$ of the computed numerical matrix ${Q^{\rho_{\ell},\delta_{\ell}}_{h_{\ell}}}$ and the computed numerical state $\mathcal{U}_{h_{\ell}} \big( Q^{\rho_{\ell},\delta_{\ell}}_{h_{\ell}} \big)$ of the algorithm  at the 500$^{\mbox{\tiny th}}$ iteration, while Figure \ref{h7} from left to right we display differences ${q^{\rho_{\ell},\delta_{\ell}}_{h_{\ell}}}_{11} - I^1_{h_{\ell}} q^\dag_{11}$, ${q^{\rho_{\ell},\delta_{\ell}}_{h_{\ell}}}_{12} - I^1_{h_{\ell}} q^\dag_{12}$, ${q^{\rho_{\ell},\delta_{\ell}}_{h_{\ell}}}_{22} - I^1_{h_{\ell}} q^\dag_{22}$ and $\mathcal{U}_{h_{\ell}} \big( Q^{\rho_{\ell},\delta_{\ell}}_{h_{\ell}} \big) - I^1_{h_\ell}\overline{u}$.

\begin{table}[H]
\begin{center}
\begin{tabular}{|c|l|l|l|l|l|l|l|}
 \hline \multicolumn{8}{|c|}{ {\bf Convergence history} }\\
 \hline
$\ell$ &\scriptsize {\bf Ite.} &\scriptsize {\bf Tol.} &\scriptsize  $\Gamma$ &\scriptsize $\Delta$ &\scriptsize $\Sigma$ &\scriptsize $\Xi$ &\scriptsize $\Lambda$\\
\hline
6& 500& 0.020795 & 9.7270e-4 & 6.4331e-4  & 0.058024  & 0.058024 & 2.1247e-4\\
\hline
12& 500& 0.005585 & 1.3063e-4  & 8.6203e-5  & 0.014679 & 0.014679 & 3.7208e-5\\
\hline
24& 500& 0.001423 & 1.6639e-5  & 1.1138e-5  & 3.6807e-3  & 3.6807e-3 & 5.4592e-6\\
\hline
48& 500& 3.5555e-4 & 2.0900e-6  & 1.4148e-6  & 9.2086e-4  & 9.2086e-4 & 7.2313e-7\\
\hline
96 & 500& 8.9001e-5 & 2.6157e-7 & 1.7825e-7  & 2.3026e-4 & 2.3026e-4  & 9.2168e-8\\
\hline
\end{tabular}
\caption{Refinement level $\ell$, number of iterations, value of tolerances, errors  $\Gamma, \Delta, \Sigma$, $\Xi$ and $\Lambda$.}
\label{b4}
\end{center}
\end{table}

\begin{table}[H]
\begin{center}
\begin{tabular}{|c|c|c|c|c|}
 \hline \multicolumn{5}{|c|}{ {\bf Experimental order of convergence} }\\
 \hline
$\ell$ &\scriptsize {\bf EOC$_\Gamma$} &\scriptsize {\bf EOC$_\Delta$} &\scriptsize {\bf EOC$_\Sigma$ = EOC$_\Xi$} &\scriptsize {\bf EOC$_\Lambda$}\\
\hline
6 & --& -- & -- &--\\
\hline
12 & 2.8965 & 2.8997 & 1.9829  & 2.5136\\
\hline
24 & 2.9728 & 2.9522 & 1.9957  & 2.7689\\
\hline
48 & 2.9930 & 2.9768 & 1.9989  & 2.9164\\
\hline
96 & 2.9982 & 2.9886 & 1.9997  & 2.9719\\
\hline
Mean of EOC & 2.9651 & 2.9543 & 1.9943 & 2.7927\\
\hline
\end{tabular}
\caption{Experimental order of convergence for $\Gamma$, $\Delta$, $\Sigma$, $\Xi$ and $\Lambda$.}
\label{b5}
\end{center}
\end{table}

\begin{figure}[H]
\begin{center}
\includegraphics[scale=0.2]{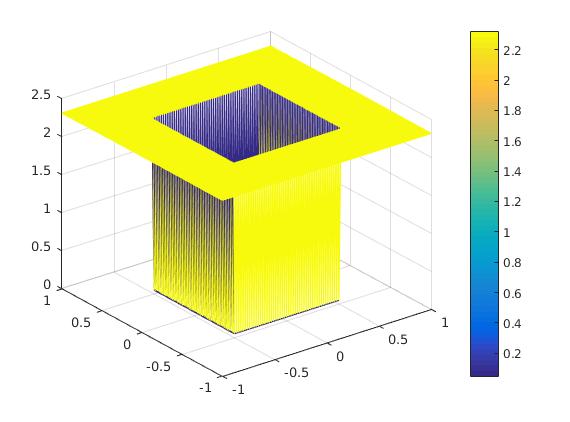}
\includegraphics[scale=0.2]{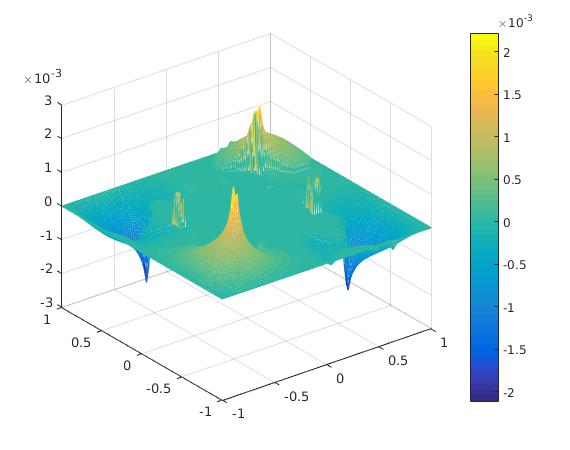}
\includegraphics[scale=0.2]{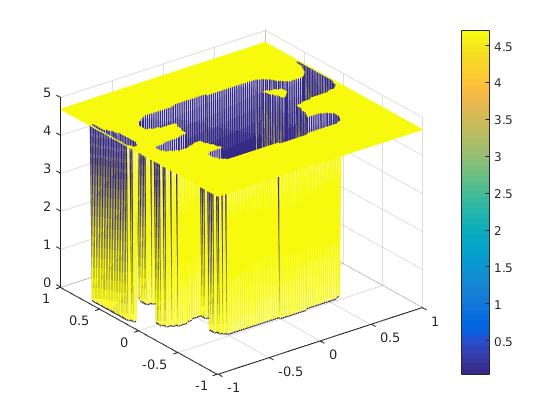}
\includegraphics[scale=0.2]{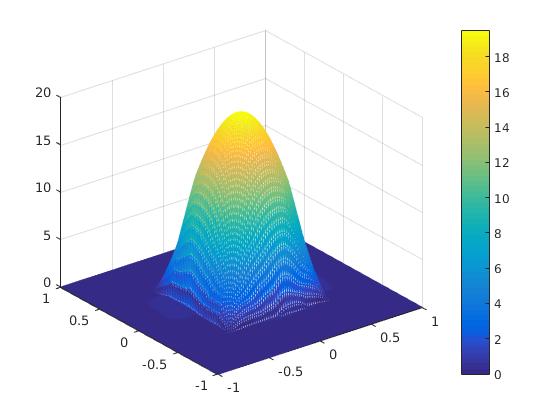}
\end{center}
\caption{Graphs of ${q^{\rho_{\ell},\delta_{\ell}}_{h_{\ell}}}_{11}, {q^{\rho_{\ell},\delta_{\ell}}_{h_{\ell}}}_{12}$, ${q^{\rho_{\ell},\delta_{\ell}}_{h_{\ell}}}_{22}$ and $\mathcal{U}_{h_{\ell}} \big( Q^{\rho_{\ell},\delta_{\ell}}_{h_{\ell}} \big)$ at the 50$^{\mbox{\tiny th}}$ iteration.}
\label{h8}
\end{figure}

Finally, Figure \ref{h8} from left to right we perform graphs of ${q^{\rho_{\ell},\delta_{\ell}}_{h_{\ell}}}_{11}, {q^{\rho_{\ell},\delta_{\ell}}_{h_{\ell}}}_{12}$, ${q^{\rho_{\ell},\delta_{\ell}}_{h_{\ell}}}_{22}$ and $\mathcal{U}_{h_{\ell}} \big( Q^{\rho_{\ell},\delta_{\ell}}_{h_{\ell}} \big)$ at the 50$^{\mbox{\tiny th}}$ iteration. At this iteration the value of tolerance is 4.1052, while errors  $\Gamma, \Delta, \Sigma, \Xi$ and $\Lambda$ are 6.9093, 3.9500, 9.9183, 46.2302 and 45.1537, respectively.

We close this section by noting that the proposed method may be extended to the case where the observation $z^\delta$ is only available in a compact subset $\Omega_{\mbox{obs}}$ of the domain $\Omega$, i.e. $\Omega_{\mbox{obs}} \Subset \Omega$. We then use a suitable $H^1_0(\Omega)$-extension $\widehat{z}^\delta$ of $z^\delta$ as measurement in our cost functional. We then consider the following strictly convex minimization problem: 
$$\min_{Q \in \mathcal{Q}_{ad}} \int_{\Omega} Q \nabla \big( \mathcal{U}(Q)
- \widehat{z}^\delta\big) \cdot \nabla \big(\mathcal{U}(Q) - \widehat{z}^\delta \big) + \rho \| Q \|^2_{{L^2(\Omega)}^{d\times d}}
\eqno \left(\widehat{\mathcal{P}}^{\rho,\delta} \right)$$
instead of $\left(\mathcal{P}^{\rho,\delta} \right)$. This problem then attains a unique solution $\widehat{Q}^{\rho,\delta}$, as in the case with full observations.

\section*{Acknowledgments} 

We thank the three referees and the associate editor for their valuable comments and suggestions. The author M.\ H.\  was supported by {\it Lothar Collatz Center for Computing in Science}. The author T.\ N.\ T.\ Q.\ was supported by {\it Alexander von Humboldt-Foundation}, Germany.


\begin{thebibliography}{10}

\bibitem{attouch} H.\ Attouch, G.\ Buttazzo and G.\ Michaille,
{\it Variational Analysis in Sobolev and BV Space}, SIAM
Mathematical Programming Society Philadelphia, 2006.

\bibitem{Bernardi1} C.\ Bernardi, Optimal finite element interpolation on curved domain, {\it SIAM J.\ Numer. Anal.} 26(1989), 1212--1240.

\bibitem{Bernardi2} C.\ Bernardi and V.\ Girault, A local regularization operator for triangular and quadrilateral finite elements, {\it SIAM J.\ Numer. Anal.} 35(1998), 1893--1916.

\bibitem{Brenner_Scott} S.\ Brenner and R.\ Scott, {\it The Mathematical Theory of Finite Element Methods}, Springer-Verlag, New York, 1994.

\bibitem{ChanTai2003}
T.\ F.\ Chan and X.\ C.\ Tai, Identification of discontinuous
coefficients in elliptic problems using total variation
regularization, {\em SIAM J.\ Sci.\ Comput.} 25(2003), 881--904.

\bibitem{Chavent_Kunisch2002}
G.\ Chavent and K.\ Kunisch, The output least squares
identifiability of the diffusion coefficient from an
$H^1$-observation in a 2-D elliptic equation, {\em ESAIM Control
Optim.\ Calc.\ Var.} 8(2002), 423--440.

\bibitem{ChenZou}
Z.\ Chen and J.\ Zou, An augmented Lagrangian method for identifying discontinuous parameters in elliptic systems, {\em SIAM J. Control Optim.} 37(1999), 892--910.

\bibitem{Cherlenyak} I.\ Cherlenyak, Numerische L\"osungen inverser
Probleme bei elliptischen Differentialgleichungen. {\it Dr.\ rer.\
nat.\ Dissertation}, Universit\"at Siegen, 2009, Verlag Dr.\ Hut,
M\"unchen 2010.

\bibitem{Chicone} C.\ Chicone and J.\ Gerlach, A note on the identifiability
of distributed parameters in elliptic equations, {\it SIAM J.\
Math.\ Anal.} 18(1987), 1378--1384.

\bibitem{Ciarlet} P.\ G.\ Ciarlet, {\it Basis Error Estimates for Elliptic Problems,} Handbook of Numerical Analisis, Vol.\ II, P.\ G.\
Ciarlet and J. -L. Lions, eds., North-Holland, Amsterdam, 1991.

\bibitem{Clement} P.\ Cl\'ement, Approximation by finite element functions
using local regularization, {\it RAIRO  Anal.\ Num\'er.} 9(1975),
77--84.

\bibitem{Deckelnick_Hinze_2011} 
K.\ Deckelnick and M.\ Hinze, Identification of matrix parameters in elliptic PDEs, {\it Control \& Cybernetics} 40(2011), 957--970.

\bibitem{Deckelnick}
K.\ Deckelnick and M.\ Hinze, Convergence and error analysis of a
numerical method for the identification of matrix parameters in
elliptic PDEs, {\em Inverse Problems} 28(2012), 115015 (15pp).

\bibitem{Engl_Hanke_Neubauer}
H.\ W.\ Engl, M.\ Hanke and A.\ Neubauer, {\em Regularization of
Inverse Problems},  Mathematics and its Applications, 375. Kluwer
Academic Publishers Group, Dordrecht, 1996.

\bibitem{EnglKuNe}
H.\ W.\ Engl, K.\ Kunisch and A.\ Neubauer, Convergence rates for
Tikhonov regularization of nonlinear ill-posed problems, {\em
Inverse Problems} 5(1989), 523--540.

\bibitem{EnglZou}
H.\ W.\ Engl and J.\ Zou, A new approach to convergence rate analysis of Tikhonov regularization for parameter identification in heat conduction, {\em Inverse Problems} 16(2000), 1907--1923.

\bibitem{Enyi}
C.\ D.\ Enyi and M.\ E.\ Soh, Modified gradient-projection algorithm
for solving convex minimization problem in Hilbert spaces,
{\em International Journal of Applied Mathematics}  44(2014), 144--150.

\bibitem{Falk}
R.\ S.\ Falk, Error estimates for the numerical identification of a
variable coefficient, {\em Math.\ Comp.} 40(1983), 537--546.

\bibitem{Flemming-Hofmann} J.\ Flemming and B.\ Hofmann, Convergence rates in constrained Tikhonov regularization: equivalence of projected source
conditions and variational inequalities, {\em Inverse Problems} 27(2011), 085001
(11pp).

\bibitem{Grasmair}
M.\ Grasmair, Generalized Bregman distances and convergence rates
for non-convex regularization methods, {\em Inverse problems} 26(2010), 115014 (16pp).

\bibitem{Haoq}
D.\ N.\ H{\`a}o and T.\ N.\ T.\ Quyen, Convergence rates for
Tikhonov regularization of coefficient identification problems in
Laplace-type equations, {\em Inverse Problems} 26(2010), 125014
(23pp).

\bibitem{hao_quyen1}
D.\ N.\ H{\`a}o and T.\ N.\ T.\ Quyen, Convergence rates for total
variation regularization of coefficient identification problems in
elliptic equations I, {\it Inverse Problems} 27(2011), 075008
(28pp).

\bibitem{hao_quyen2}
D.\ N.\ H{\`a}o and T.\ N.\ T.\ Quyen, Convergence rates for total
variation regularization of coefficient identification problems in
elliptic equations II, {\it J.\ Math.\ Anal.\ Appl.} 388(2012),
593--616.

\bibitem{hao_quyen3}
D.\ N.\ H{\`a}o and T.\ N.\ T.\ Quyen, Convergence rates for
Tikhonov regularization of a two-coefficient identification problem
in an elliptic boundary value problem, {\em Numer.\ Math.}
120(2012), 45--77.

\bibitem{hao_quyen4}
D.\ N.\ H{\`a}o and T.\ N.\ T.\ Quyen, Finite element methods for
coefficient identification in an elliptic equation, {\em Appl.\
Anal.} 93(2014), 1533--1566.

\bibitem{Hinze}
M.\ Hinze, A variational discretization concept in control constrained optimization: the linear-quadratic case, {\em Comput. Optim. Appl.} 30(2005), 45--61.

\bibitem{Hoffmann_Sprekels}
K.\ H.\ Hoffmann and J.\ Sprekels, On the identification of
coefficients of elliptic problems by asymptotic regularization, {\em
Numer.\ Funct.\ Anal.\ Optim.} 7(1985), 157--177.

\bibitem{HKPS}
B.\ Hofmann, B.\ Kaltenbacher C.\ P\"oschl and O Scherzer, A convergence rates result for Tikhonov regularization in Banach spaces with non-smooth operators, {\em Inverse Problems} 23(2007), 987--1010.

\bibitem{Hohage}
T.\ Hohage and F.\ Weidling, Verification of a variational source condition for acoustic inverse medium scattering problems,
{\em Inverse Problems} 31(2015), 075006 (14pp).

\bibitem{HohageWerner}
T.\ Hohage and F.\ Werner, Convergence Rates for Inverse Problems with Impulsive Noise, {\em SIAM J.\ Numer.\ Anal.} 52(2014), 1203--1221.

\bibitem{Kaltenbacher_Schoberl}
B.\ Kaltenbacher and J.\ Sch\"oberl, A saddle point variational
formulation for projection-regularized parameter identification,
{\em Numer. Math.} 91(2002), 675--697.

\bibitem{KeungZou98}
Y.\ L.\ Keung and J.\ Zou, Numerical identifications of parameters in parabolic systems, {\em Inverse Problems} 14(1998), 83--100.

\bibitem{KeungZoz2000}
Y.\ L.\ Keung and J.\ Zou, An efficient linear solver for nonlinear parameter identification problems, {\em SIAM J. Sci. Comput.} 22(2000), 1511--1526.

\bibitem{Know2}
I.\ Knowles, Uniqueness for an elliptic inverse problem, {\em SIAM
J.\ Appl.\ Math.} 59(1999), 1356--1370.

\bibitem{kolo}
R.\ V.\ Kohn and B.\ D.\ Lowe, A variational method for parameter
identification, {\em RAIRO Mod\'el Math.\ Anal.\ Num\'er.} 22(1988),
119--158.

\bibitem{Kohn_Vogelius1}
R.\ V.\ Kohn and M.\ Vogelius, Determining Conductivity by Boundary
Measurements, {\em Comm.\ Pure Appl.\ Math.} 37(1984), 289--298.

\bibitem{Kohn_Vogelius2}
R.\ V.\ Kohn and M.\ Vogelius, Relaxation of a variational method for impedance computed tomography, {\em Comm.\ Pure Appl.\ Math.}
40(1987), 745--777.

\bibitem{Kunisch}
K.\ Kunisch, {\it Numerical methods for parameter estimation problems}, Inverse problems in diffusion processes (Lake St. Wolfgang, 1994), 199--216, SIAM, Philadelphia, PA, 1995.

\bibitem{Murat_Tartar} F.\ Murat and L.\ Tartar, {\it H-Convergence.}
Topics in the Mathematical Modelling of composite materials,
Andrej Charkaev and Robert Kohn, eds., pp. 21--43, Birkh{\"a}user,
1997.

\bibitem{Rannacher_Vexler}
R.\ Rannacher and B.\ Vexler, A priori error estimates for the
finite element discretization of elliptic parameter identification
problems with pointwise measurements, {\em SIAM J.\ Control Optim.}
44(2005), 1844-1863.

\bibitem{Ric}
G.\ R.\ Richter, An inverse problem for the steady state diffusion
equation {\em SIAM J.\ Appl.\ Math.}, 41(1981), 210--221.

\bibitem{scott_zhang} R.\ Scott and S.\ Y.\ Zhang,  Finite element
interpolation of nonsmooth function satisfying boundary conditions,
{\it Math.\ Comp.} 54(1990), 483--493.

\bibitem{Spagnolo} S.\ Spagnolo, Sulla convergenza di
soluzioni di equazioni paraboliche ed ellittiche. {\it Ann.\ Sc.\
Norm.\ Sup.\ Pisa} 22(1968), 571--597.

\bibitem{Tartar} L.\ Tartar, {\it Estimation of homegnized coefficients.}
Topics in the Mathematical Modelling of composite materials,
Andrej Charkaev and Robert Kohn, eds., pp. 9--20, Birkh{\"a}user,
1997.

\bibitem{tro}
G.\ M.\ Troianiello, {\it Elliptic Differential Equations and
Obstacle Problems}, Plenum, New York, 1987.

\bibitem{Vainikko-Kunisch}
G.\ Vainikko and K.\  Kunisch, Identifiability of the transmissivity coefficient in an elliptic boundary value problem, {\it Z. Anal. Anwendungen} 12(1993), 327--341.

\bibitem{wang_zou}
L. Wang and J. Zou, Error estimates of finite element methods for
parameter identification problems in elliptic and parabolic systems,
{\em Discrete Contin.\ Dyn.\ Syst.\ Ser.\ B.} 14(2010), 1641--1670.

\bibitem{WernerHohage2012}
F.\ Werner and T.\ Hohage, Convergence rates in expectation for Tikhonov-type
regularization of inverse problems with Poisson data, {\em Inverse problems} 28(2012), 104004 (15pp).

\end{thebibliography}
\end{document}